\documentclass[a4paper, 12pt]{amsart}
  \usepackage{latexsym}
  \usepackage[all]{xy}
  \usepackage{amsfonts}
  \usepackage{amsthm}
  \usepackage{amsmath}
  \usepackage{amssymb}
  \usepackage[mathscr]{eucal}
  \usepackage{enumerate}
  \usepackage[latin1]{inputenc}
  \usepackage{a4wide}

\newcommand{\naturales}{\mathbb{N}}

 \newcommand{\rmod}[1]{\mathcal{M}_{#1}}

\newcommand{\rcomod}[1]{\mathcal{M}^{#1}}

\newcommand{\bimod}[2]{{}_{#1}\mathcal{M}_{#2}}

\newcommand{\ring}[1]{\mathscr{#1}}
\newcommand{\coring}[1]{\mathfrak{#1}}
\newcommand{\tensor}[1]{\otimes_{#1}}

\newcommand{\bara}[1]{\overline{#1}}

\newcommand{\cC}{\coring{C}}

\newcommand{\cD}{\coring{D}}

\newcommand{\idemp}[1]{\mathrm{Idemp}(#1)}
\newcommand{\Bk}{\mathbb{K}}
\newcommand{\lr}[1]{\langle #1 \rangle}
\renewcommand{\hom}[3]{\mathrm{Hom}_{#1}(#2,#3)}
\newcommand{\rend}[2]{\mathrm{End}(#1_{#2})}

\newcommand{\dlimit}[1]{\underset{\longrightarrow}{\mathrm{lim}}#1}
\newcommand{\cat}[1]{\mathcal{#1}}
\newcommand{\coanillos}[1]{#1-\mathsf{coring}}
\newcommand{\comonads}[1]{\cat{#1}-\mathsf{comonad}}
\newcommand{\msf}[1]{\mathsf{#1}}
\newcommand{\bd}[1]{\boldsymbol{#1}}
\newcommand{\udual}[1]{{#1}^{\dag}}
\newcommand{\fk}[1]{\mathfrak{#1}}
\newcommand{\LR}[1]{\left(\underset{}{} #1 \right)}
\newcommand{\td}[1]{\widetilde{#1}}

\xyoption{2cell} %\xycompileto{D}
  \newtheorem{proposition}{Proposition}[section]
  \newtheorem{lemma}[proposition]{Lemma}
  
  \newtheorem{corollary}[proposition]{Corollary}
  \newtheorem{theorem}[proposition]{Theorem}

  \theoremstyle{definition}
  \newtheorem{definition}[proposition]{Definition}
  \newtheorem{example}[proposition]{Example}

  \theoremstyle{remark}
  \newtheorem{remark}[proposition]{Remark}

\begin{document}

\baselineskip 16pt

\title[Corings over rings with local units]{Corings over rings with local units}
\author[L. El Kaoutit]{L. El Kaoutit} 
\address{Universidad de Granada. Departamento de \'Algebra.
Facultad de Educaci\'on y Humanidades de Ceuta. El Greco N10, 51002
Ceuta, Spain.}

\date{March 2006}


\begin{abstract}
We show that the bicategory whose $0$-cells are corings over rings
with local units is bi-equivalent to the bicategory of comonads in
(right) unital modules whose underlying functors are right exact and
preserve direct sums. A base ring extension of a coring by an
adjunction is introduced as well.
\end{abstract}

\maketitle

\section{Introduction}
Corings over rings with identity have been intensively studied in the
last years. A detailed discussion can be found in
\cite{Brzezinski/Wisbauer:2003}. It is well known that any coring
entails a comonad (i.e., cotriple) in the category of right modules
over the base ring. The converse is also true for the case of
coalgebras, taking functors which are right exact and preserve direct
sums, as was checked in \cite{Eilenberg/Moore:1965}. Thus, the
categorical study of corings and their comodules is entirely linked
to the study of certain comonads and their universal cogenerators
(see \cite{Eilenberg/Moore:1965}) in non necessarily monoidal
categories. For instance, many new examples of corings can be built
using earlier constructions in comonads theory (e.g. the
distributive laws of J. Beck \cite{Beck:1969}), or just by
considering comonads with a special cogenerator.

In this paper, we study these relationships in the context of corings
over rings with local units in the sense of
\cite{Abrams:1983,Anh/Marki:1983}. This class of rings arose
naturally in the definition of infinite comatrix corings, see
\cite{Kaoutit/Gomez:2004b}, \cite{Gomez/Vercruysse:2005} and
\cite{Caenepeel/DeGroot/Vercruysse:2005}.\\
In what follows an additive covariant functor is said to be \emph{continuous}
if it is right exact and preserves direct sums.

The paper is organized as follows. Section \ref{section0} is rather
technical, and it is devoted to introduce a $2$-category with
objects ($0$-cells) comonads in Grothendieck categories, using
earlier results from \cite{Lack/Street:2002}. In Section
\ref{sectionI}, we extend a result of Watts \cite{Watts:1960} to the
case of rings with local units (Lemma \ref{ISO}), and use this to
prove that each comonad induces a coring whenever its underlying
functor is continuous (Proposition \ref{comonad2}). In Section
\ref{section4}, we establish a bi-equivalence between the bicategory
of unital bimodules and the bicategory whose $0$-cells are rings
with local units and having the categories of continuous functors
over unital right modules as Hom-Categories (Proposition \ref{B-L},
compare with \cite[Proposition 2.1]{Grandjean/Vitale:1998}). We,
then, deduce an equivalence of categories between the category of
corings over a fixed ring with local units and the category of
comonads in right unital modules whose underlying functors are
continuous (Corollary \ref{corolario}). Lastly, we derive a
bi-equivalence between the bicategory of all corings over rings with
local units and the bicategory of comonads in right unital modules
with continuous underlying functors (Theorem \ref{F-G}). The left
version of these results is similarly obtained and will not be
considered. In Section \ref{sect.app} we apply results from Sections
\ref{section0} and \ref{sectionI} to introduce a base ring extension
of a coring by an adjunction.

Throughout $\Bk$ denotes a commutative ring with identity $1$.

\section{Review on the 2-category of comonads}\label{section0}
In this section we observe that Grothendieck categories form a class
of objects ($0$-cells) in a $2$-category whose $1$-cells are
continuous functors. These will be
needed to introduce the $2$-category of comonads using the formalism
of \cite{Lack/Street:2002}. For general definitions of bicategories
and their homomorphisms we refer the reader to the fundamental paper
\cite{Benabou:1967}.
\bigskip

Recall from \cite{Eilenberg/Moore:1965} that a comonad (or cotriple)
in a category $\cat{A}$ is a three-tuple $(F, \delta,\xi)$
consisting of a functor $F: \cat{A} \rightarrow \cat{A}$ and natural
transformations $\delta_{-}: F(-) \rightarrow F \,\circ\,
F(-)=F^2(-)$, $\xi_-: F(-) \rightarrow id_{\cat{A}}(-)$ such that
\begin{equation}\label{comonad}
\xymatrix@R=30pt@C=60pt{F(-)\ar@{->}^-{\delta_{-}}[r]
\ar@{->}_-{\delta_{-}}[d] & F^2(-)
\ar@{->}^-{\delta_{F(-)}}[d]\\
F^2(-) \ar@{->}^-{F(\delta_{-})}[r] & F^3(-)}\qquad
\xymatrix@R=30pt@C=60pt{F(-) \ar@{->}^-{\delta_{-}}[r]
\ar@{->}_-{\delta_{-}}[d] \ar@{=}[dr] & F^2(-)
\ar@{->}^-{F(\xi_{-})}[d] \\
F^2(-) \ar@{->}^-{\xi_{F(-)}}[r] & F(-) }
\end{equation}
are commutative diagrams. It is well known from
\cite{Huber:1961} that any adjunction $\xymatrix{ S:
\mathcal{B} \ar@<0,5ex>[r] & \ar@<0,5ex>[l] \cat{A}: T}$ with $S$
left adjoint to $T$ (notation $S \dashv T$) induces a comonad $(ST, S\eta_{T},
\zeta)$ in $\cat{A}$, where $\eta: id_{\cat{B}} \rightarrow TS$ and
$\zeta:ST \rightarrow id_{\cat{A}}$ are, respectively, the unit and
the counit of this adjunction. Using the terminology of
\cite[Proposition 2.1]{Eilenberg/Moore:1965}, we say that the adjunction $S \dashv
T$ cogenerates the comonad $(ST, S\eta_{T}, \zeta)$. Now, the dual
version of \cite[Theorem 2.2]{Eilenberg/Moore:1965} asserts that for
any comonad $(F,\delta,\xi)$ in $\cat{A}$, there exists a universal
cogenerator, that is, a category $\cat{A}^F$ and an adjunction
$\xymatrix{ S^F: \cat{A}^F \ar@<0,5ex>[r] & \ar@<0,5ex>[l] \cat{A}:
T^F}$ cogenerating $(F,\delta,\xi)$ with the following universal
property: If $\xymatrix{ S: \mathcal{B} \ar@<0,5ex>[r] &
\ar@<0,5ex>[l] \cat{A}: T}$ is another cogenerator of
$(F,\delta,\xi)$, then there exists a unique functor $L: \cat{B}
\rightarrow \cat{A}^F$ such that $S^F \,\circ\, L = S$ and $L
\,\circ\, T = T^F$. The objects of $\cat{A}^F$ are referred to as
\emph{comodules}, they are pairs $(X,d^X)$ consisting of an object
$X \in \cat{A}$ and a morphism $d^X: X \rightarrow F(X)$ in
$\cat{A}$ such that
\begin{equation}\label{coacion}
\xymatrix@R=30pt@C=60pt{ X \ar@{->}^-{d^X}[r] \ar@{->}_-{d^X}[d] &
F(X) \ar@{->}^-{\delta_X}[d]
\\ F(X) \ar@{->}^-{F(d^X)}[r] & F^2(X) } \quad \xymatrix@R=30pt@C=60pt{
X \ar@{->}^-{d^X}[r] \ar@{=}[dr] & F(X) \ar@{->}^-{\xi_X}[d] \\ & X
}
\end{equation}
are commutative diagrams. A morphism in $\cat{A}^F$ is a morphism
$f: X \rightarrow X'$ in $\cat{A}$, where $(X,d^X)$, $(X',d^{X'}) \in
\cat{A}^F$ such that
\begin{equation}\label{colineal}
\xymatrix@R=30pt@C=60pt{X \ar@{->}_-{d^X}[d] \ar@{->}^-{f}[r] &
X' \ar@{->}^-{d^{X'}}[d] \\
F(X) \ar@{->}^-{F(f)}[r] & F(X') }
\end{equation}
is a commutative diagram. The functor $S^F$ is then the forgetful
functor $S^F(X,d^X)=X$ for every comodule $(X,d^X)$, and the functor
$T^F$ is defined by $T^F(Y) = (F(Y), \delta_Y)$, $T^F(f)=F(f)$ for
any object $Y$ and morphism $f$ of $\cat{A}$.

Given $(F,\delta,\xi)$ a comonad in $\cat{A}$, the Kleisli
\cite{Kleisli:1965} category $\cat{A}_{K}$, is defined as follows.
The objects of $\cat{A}_K$ are the same as those of $\cat{A}$. For
each pair of objects $Y,Y'$, the morphism set is defined by
$$\hom{\cat{A}_K}{Y}{Y'} \,\, := \, \, \hom{\cat{A}}{F(Y)}{Y'}.$$
The composition law and identities are canonically derived using
$\delta$ and $\xi$. As was proved in \cite{Kleisli:1965}, there is an
adjunction $\xymatrix{ S: \cat{A}_K \ar@<0,5ex>[r] & \ar@<0,5ex>[l]
\cat{A}: T}$ cogenerating $(F,\delta,\xi)$, where $S(Y)=F(Y)$ and
$S(f)= F(f) \,\circ\, \delta_Y$ for any object $Y$ and morphism
$f$ of $\cat{A}_K$, and $T(X)=X$, $T(g)= \xi_{X'} \,\circ\, F(g)$
for any objects $X$, $X'$ and morphism $g:X\rightarrow X'$ of
$\cat{A}$. The unique factorizing functor is given as follows. The
adjunction $S^F \dashv T^F$ gives a natural isomorphism
$$ \hom{\cat{A}_K}{X}{X'} \,\, \cong \,\,
\hom{\cat{A}^F}{T^F(X)}{T^F(X')},$$ which defines a functor $L:
\cat{A}_K \rightarrow \cat{A}^F$ by $L(X) = T^F(X)=
(F(X),\delta_X)$, for any object $X \in \cat{A}_K$, and $L(f)=
F(f) \,\circ\, \delta_X: F(X) \rightarrow F(X')$, for any morphism
$f \in \hom{\cat{A}_K}{X}{X'}$. Clearly $L \,\circ\, T=T^F$, $S^F
\,\circ\, L =S$, and $L$ is a fully faithful functor.

\begin{lemma}\label{dos}
Consider an adjunction $\xymatrix{S: \cat{A} \ar@<0,5ex>[r] &
\ar@<0,5ex>[l] \cat{B}:T}$, $S\dashv T$ with unit $\eta:
id_{\cat{A}} \rightarrow TS$ and counit $\zeta: ST \rightarrow
id_{\cat{B}}$ together with $(F,\delta,\xi)$ a comonad in $\cat{A}$.
Then
\begin{enumerate}
\item $[S\dashv T](F,\delta,\xi)\,\, :=\,\,(SFT, S F \eta_{FT}
\,\circ\, S\delta_{T}, \zeta \,\circ\, S \xi_T)$ is a comonad in
$\cat{B}$.

\item There is a functor $\ring{S}: \cat{A}^F \longrightarrow
\cat{B}^{SFT}$ defined by
$$
(X,d^X) \rightarrow \left(\underset{}{}S(X),d^{S(X)}= SF\eta_X
\,\circ\, Sd^X \right), \qquad (f \rightarrow S(f)).$$
\end{enumerate}
\end{lemma}
\begin{proof}
(1) This is \cite[Theorem 4.2]{Huber:1961}.\\ (2) It is clear that
we have an adjunction $$\xymatrix@C=60pt{ \cat{A}^F \ar@<0,5ex>^-{S
\,\circ\, S^F}[r] & \ar@<0,5ex>^-{T^F \,\circ\, T}[l] \cat{B} }$$
with $S \,\circ\, S^F \dashv T^F \,\circ\, T$ which cogenerates the
comonad $(SFT, S F \eta_{FT} \,\circ\, S\delta_{T}, \zeta \,\circ\,
S \xi_T)$. Following the proof of the dual version of \cite[Theorem
2.2]{Eilenberg/Moore:1965}, we obtain that the stated functor is the
well known unique factorization functor.
\end{proof}

Following \cite{Bautista/Colavita/Salmeron:1981}, a morphism between
two comonads $(F,\delta,\xi)$ and $(F',\delta',\xi')$ in a category
$\cat{A}$ is a natural transformation $\Phi_{-}: F(-) \rightarrow
F'(-)$ which turns commutative the following diagrams
\begin{equation}\label{compatible}
\xymatrix@R=30pt@C=60pt{ F(-) \ar@{->}^-{\Phi_{-}}[r]
\ar@{->}_-{\delta_{-}}[d]
 & F'(-) \ar@{->}^-{\delta'_{-}}[d] \\
F^2(-) \ar@{->}^-{\varrho(\Phi)_{-}}[r] & (F')^2(-) } \qquad
\xymatrix@R=30pt@C=60pt{ F(-) \ar@{->}^-{\Phi_{-}}[r]
\ar@{->}_-{\xi_{-}}[d] & F'(-) \ar@{->}^-{\xi'_{-}}[d] \\
id_{\cat{A}}(-) \ar@{=}[r] & id_{\cat{A}}(-), }
\end{equation}
where $\varrho(\Phi)_{-}$ is the natural transformation defined by
\begin{equation}\label{varrho}
\varrho(\Phi)_{-} \,\, := \,\, F'\Phi_{-} \, \,\circ\, \,
\Phi_{F(-)} \,\, = \,\, \Phi_{F'(-)} \, \,\circ\, \, F\Phi_{-}.
\end{equation}
Given a morphism of comonads $\Phi: F \rightarrow F'$, we have an
induction functor $$ \xymatrix{ (-)_{\Phi}: \cat{A}^F \ar@{->}[r] &
\cat{A}^{F'} }$$ defined by $(X,d^X)_{\Phi} = (X, \Phi_X \,\circ\,
d^X)$ for any $(X,d^X) \in \cat{A}^F$ and the identity on
morphisms.

\begin{example}\label{ejemplo}
Let $(F,\delta,\xi)$ be any comonad in a category $\cat{A}$. Clearly
$ \xi: F \rightarrow id_{\cat{A}}$ is morphism of comonads, where
$id_{\cat{A}}$ is endowed with a trivial comonad structure. If we
consider any adjunction $\xymatrix{S: \cat{A} \ar@<0,5ex>[r] &
\ar@<0,5ex>[l] \cat{B}:T}$  with $S \dashv T$ and $[S \dashv
T](F,\delta,\xi)$ the associated comonad of Lemma \ref{dos}, then
we can easily check that $S\xi_T: SFT \rightarrow ST$ is also a
morphism of comonads in $\cat{B}$, where $(ST,S\eta_T,\zeta)$ is the
comonad cogenerated by $S \dashv T$.
\end{example}

Next, we are going to look at the case in which certain comonads and their morphisms
form a set-category and try to interpret the above constructions by
means of functors between those categories. The following definition
makes sense after \cite[Lemma 5.1]{Iglesias/Gomez/Nastasescu:1999},
where it was proved that, over Grothendieck category, natural
transformations between continuous functors form a set.

\begin{definition}\label{continuos}
Let $\cat{A}$ be a Grothendieck category. We define the category of
comonads in $\cat{A}$ and denote it by $\comonads{A}$, as the
category whose objects are comonads $(F,\delta,\xi)$ in $\cat{A}$
with $F: \cat{A} \rightarrow \cat{A}$ a continuous functor (thus it
also preserves inductive limits); and whose morphisms are natural
transformations satisfying the commutativity of the diagrams
\eqref{compatible}.
\end{definition}

With this definition we can give an elegant restatement of Lemma
\ref{dos}(1).

\begin{proposition}\label{funtor}
Let $\cat{A}$ and $\cat{B}$ be two Grothendieck categories together
with an adjunction $\xymatrix{S: \cat{A} \ar@<0,5ex>[r] &
\ar@<0,5ex>[l] \cat{B}:T}$, $S\dashv T$ whose unit and counit are,
respectively, $\eta: id_{\cat{A}} \rightarrow TS$ and $\zeta: ST
\rightarrow id_{\cat{B}}$. Assume that $S$ and $T$ are continuous functors.
\begin{enumerate}
\item The assignment of Lemma \ref{dos}(1)
$$ (F,\delta,\xi) \rightarrow [S\dashv T](F,\delta,\xi),
\quad [\Phi: F \rightarrow F'] \rightarrow \left[\underset{}{}
[S\dashv T](\Phi)\,=\,[S\Phi_T: SFT \rightarrow SF'T]\right] $$
defines a functor $[S\dashv T]: \comonads{A} \longrightarrow
\comonads{B}$.

\item Given $\xymatrix{P: \cat{B} \ar@<0,5ex>[r] & \ar@<0,5ex>[l]
\cat{C}:Q}$ another adjunction where $\cat{C}$ is a Grothendieck
category and $P$, $Q$ are continuous functors. Then, we have the following composition
$$ [P\dashv Q] \, \,\circ\, \, [S\dashv T] \,\,=\,\, [PS \dashv
TQ],$$ where $[P\dashv Q]: \comonads{B} \rightarrow \comonads{C}$
and $[PS\dashv TQ]: \comonads{A} \rightarrow \comonads{C}$ are
functors defined as above.
\end{enumerate}
\end{proposition}
\begin{proof}
(1) We only show that $$S\Phi_{T}=\bara{\Phi}: \,\,\bara{F}=SFT
\longrightarrow \bara{F'}=SF'T$$ is a morphism of comonads where
$(F,\delta,\xi)$ and $(F',\delta',\xi')$ are comonads in $\cat{A}$.
We have
$$  \zeta \,\circ\, S\xi'_T \,\circ\, \bara{\Phi} \,\,=\,\, \zeta \,\circ\, S\xi'_T \,\circ\, S\Phi_T
\,\, = \,\, \zeta \,\circ\, S \left( \underset{}{} \xi' \,\circ\,
\Phi \right)_T  \,\, =\,\, \zeta \,\circ\, S\xi_T$$ which shows that
the second diagram in equation \eqref{compatible} commutes. On the
other hand, we have

\begin{eqnarray*}
  SF'\eta_{F'T} \,\circ\, S\delta'_T \,\circ\, S\Phi_T &=&
  S \left( \underset{}{} F'\eta_{F'} \,\circ\, \delta' \,\circ\, \Phi \right)_T \\
   &=& S \left( \underset{}{} F'\eta_{F'} \,\circ\, \varrho(\Phi) \,\circ\, \delta
   \right)_T,\quad \text{ by }\eqref{compatible} \\
   &=& S \left( \underset{}{} F'\eta_{F'} \,\circ\, F'\Phi \,\circ\, \Phi_F \,\circ\, \delta
   \right)_T ,\quad \text{ by }\eqref{varrho}  \\
   &=& S \left( \underset{}{} F'(\eta_{F'} \,\circ\, \Phi) \,\circ\, \Phi_F \,\circ\, \delta
   \right)_T \\
   &=& S \left( \underset{}{} F'TS\Phi \,\circ\, F'\eta_F \,\circ\, \Phi_F \,\circ\, \delta
   \right)_T,\quad \eta_{-} \text{ is natural} \\
   &=& S \left( \underset{}{} F'TS\Phi \,\circ\, \Phi_{TSF} \,\circ\, F\eta_F \,\circ\, \delta
   \right)_T,\quad \Phi_{-} \text{ is natural} \\
   &=& SF'TS\Phi_T \,\circ\, S\Phi_{TSFT} \,\circ\, SF\eta_{FT} \,\circ\, S\delta_T \\
   &=& \varrho(S\Phi_T) \,\circ\, SF\eta_{FT} \,\circ\, S\delta_T.
\end{eqnarray*}
Thus the first diagram in equation \eqref{compatible} commutes for
$\bara{F}$, $\bara{F'}$ and $\bara{\Phi}=S\Phi_T$.\\
(2) Straightforward.
\end{proof}

Using Definition \ref{continuos}, one can adapt the formalism
introduced in \cite{Lack/Street:2002} (see also \cite{Street:1972})
for monads in arbitrary $2$-category, to the setting of Grothendieck
categories as follows: First, using \cite[Lemma
5.1]{Iglesias/Gomez/Nastasescu:1999}, we obtain a $2$-category
constructed by the following data:
\begin{enumerate}[$\bullet$]
\item \emph{Objects ($0$-cells)}: All Grothendieck categories.

\item \emph{$1$-cells}: An $1$-cell from $\cat{B}$ to $\cat{A}$ is a
continuous functor $F: \cat{A} \rightarrow \cat{B}$.

\item \emph{$2$-cells}: Natural transformations.
\end{enumerate}
Associated to this $2$-category, we construct, as in
\cite{Lack/Street:2002}, the right Eilenberg-Moore $2$-category of
comonads:

\begin{enumerate}[$\bullet$]

\item \emph{Objects} (\emph{$0$-cells}): They are pairs
$(F,\delta,\xi: \cat{A})$ consisting of a Grothendieck category
$\cat{A}$ and a comonad $\bd{F}=(F,\delta,\xi)$ in $\cat{A}$ such
that $F: \cat{A} \rightarrow \cat{A}$ is a continuous
functor (i.e., $\bd{F}$ is an object of
the category $\comonads{A}$ of Definition \ref{continuos}).

\item \emph{$1$-cells}: An $1$-cell from $(\bd{G}:\cat{B})$ to
$(\bd{F}:\cat{A})$ (here $\bd{G}=(G,\Omega,\gamma) \in
\comonads{B}$), is a pair $(S,\msf{s})$ consisting of a continuous
functor $S:\cat{A} \rightarrow \cat{B}$, and a natural transformation $\msf{s}: SF
\longrightarrow GS$ satisfying the commutativity of the following
two diagrams
\begin{equation}\label{msf}
\xymatrix@R=30pt@C=40pt{
SF\ar@{->}^-{\msf{s}}[r] \ar@{->}|-{S\xi}[d] & GS \ar@{->}|-{\gamma_S}[d] \\
S \ar@{=}[r] & S,} \qquad \xymatrix@R=30pt@C=60pt{SF
\ar@{->}^-{\msf{s}}[r] \ar@{->}|-{S\delta}[d] & GS \ar@{->}^-{\Omega_S}[r] & G^2S \\
SF^2 \ar@{->}^-{\msf{s}_F}[r] & GSF \ar@{->}|-{G\msf{s}}[ur] & }
\end{equation}The identity $1$-cell for a given object
$(\bd{F}:\cat{A})$ is provided by $(id_{\cat{A}}, id_{F(-)})$.

\item \emph{$2$-cells}: Given $(S,\msf{s})$ and $(S',\msf{s}')$
two $1$-cells from $(\bd{G}:\cat{B})$ to $(\bd{F}:\cat{A})$, a
$2$-cell $(S,\msf{s}) \rightarrow (S',\msf{s}')$ is a natural
transformation $\alpha: SF \rightarrow S'$ turning commutative the
following diagram
\begin{equation}\label{dos-cells}
\xymatrix@R=30pt@C=60pt{ SF \ar@{->}^-{S\delta}[r]
\ar@{->}|-{S\delta}[d] & SF^2 \ar@{->}^-{\msf{s}_F}[r] & GSF \ar@{->}|-{G\alpha}[d] \\
SF^2 \ar@{->}^-{\alpha_F}[r] & S'F \ar@{->}^-{\msf{s}'}[r] & GS'. }
\end{equation}
\end{enumerate}
The category constructed by all $1$ and $2$-cells from
$(\bd{G}:\cat{B})$ to $(\bd{F}:\cat{A})$ will be denoted by
${}_{\bd{F}}\ring{C}_{\bd{G}}$. The laws composition are given as
follows. Let $(S,\msf{s})$, $(S',\msf{s}')$, and $(S'',\msf{s}'')$
be three $1$-cells from $(\bd{G}:\cat{B})$ to $(\bd{F}:\cat{A})$
with $2$-cells $\alpha: (S,\msf{s}) \rightarrow (S',\msf{s}')$ and
$\alpha': (S',,\msf{s}') \rightarrow (S'',,\msf{s}'')$. Then
\begin{equation}\label{comp-C}
\alpha' \,\,\underline{\circ}\, \, \alpha \,\,=\,\, \alpha'\,
\,\circ\, \, \alpha_F \, \,\circ\, S{\delta},\quad \text{ where
}\bd{F}=(F,\delta,\xi).
\end{equation}
Given $(P,\msf{p})$ and $(P',\msf{p}')$ two $1$-cells from
$(\bd{H}:\cat{C})$ to $(\bd{G}:\cat{B})$, together with $2$-cells
$\alpha:(S,\msf{s}) \rightarrow (S',\msf{s}')$ and
$\beta:(P,\msf{p}) \rightarrow (P',\msf{p}')$, the vertical
composition is given by
\begin{equation}\label{vertical}
(S,\msf{s}) \, .\, (P,\msf{p}) \,\, = \,\, (PS, \msf{p}_S \,\circ\,
P\msf{s}), \text{ and } (S',\msf{s}') \,.\, (P',\msf{p}') \,\, =
\,\, (P'S', \msf{p}'_{S'} \,\circ\, P'\msf{s}')
\end{equation}
The horizontal composition $  \alpha \,.\, \beta: (PS,
\msf{p}_S \,\circ\, P\msf{s}) \rightarrow (P'S', \msf{p}'_{S'}
\,\circ\, P'\msf{s}')$ is defined by
\begin{equation}\label{horizontal}
\xymatrix@R=30pt@C=60pt{ PSF \ar@{->}^-{PS\delta}[r]
\ar@{-->}|-{\alpha\,.\, \beta}[drrr]
& PSF^2 \ar@{->}^-{P\msf{s}_F}[r] & PGSF \ar@{->}^-{PG\alpha}[r] & PGS' \ar@{->}|-{\beta_{S'}}[d] \\
& & & P'S'}
\end{equation}

Associated to an $1$-cell $(S,\msf{s}) \in
{}_{\bd{F}}\ring{C}_{\bd{G}}$, there is a functor connecting the
universal cogenerators. Namely, there is an additive functor
\begin{equation}\label{Snegra}
\xymatrix@R=0pt@C=60pt{\bd{\ring{S}}: \, \cat{A}^F  \ar@{->}[r]&
\cat{B}^G, }
\end{equation}
sending
$$\left( (X,d^X) \longrightarrow (S(X),\msf{s}_X\,\circ\, Sd^X)
\right),\quad \left( f \rightarrow S(f) \right),$$ which clearly
turns commutative the following diagram
\begin{equation}\label{Snegra1}
\xymatrix@C=40pt{ \cat{A}^F \ar@{->}^-{\bd{\ring{S}}}[r]
\ar@{->}_-{S^F}[d] & \cat{B}^G \ar@{->}^-{S^G}[d]
\\ \cat{A} \ar@{->}^-{S}[r] & \cat{B} }
\end{equation}

As in the case of an arbitrary $2$-category \cite{Lack/Street:2002},
one can substitute the above $2$-cells (reduced forms) by the
unreduced forms, that is,  natural transformations of the form
$\alpha: SF \rightarrow GS'$ satisfying adequate conditions. The
bijection between reduced forms and unreduced forms established in
\cite{Lack/Street:2002} for monads in $2$--category, is interpreted
in our setting by the forthcoming proposition whose proof is based
upon the following. Recall that an object $V$ of an additive
category $\cat{G}$ with direct sums and cokernels, is said to be a
\emph{subgenerator}, if every object of $\cat{G}$ is a sub-object of
a $V$-generated one.
\begin{lemma}\label{additiva}
Let $\cat{A}$ be a Grothendieck category and $(F,\delta,\xi)$ be an
object of the category $\comonads{A}$. Then $\cat{A}^F$ is an
additive category with direct sums and cokernels. Furthermore, if
$U$ is a generator of $\cat{A}$ then $(F(U),\delta_U)$ is a
subgenerator of $\cat{A}^F$.
\end{lemma}
\begin{proof}
It is immediate since $F$ preserves direct sums and cokernels.
\end{proof}

\begin{proposition}\label{nat-set}
Let $\cat{A}$ and $\cat{B}$ be two Grothendieck categories, and
$\bd{F}=(F,\delta,\xi) \in \comonads{A}$, $\bd{G}=(G,\Omega,\gamma)
\in \comonads{B}$. Considering $(S,\msf{s})$ and $(S',\msf{s}')$ two
$2$-cells from $(\bd{G}:\cat{B})$ to $(\bd{F}:\cat{A})$, with the
associated functors $\bd{\ring{S}}, \bd{\ring{S}'}: \cat{A}^F
\rightarrow \cat{B}^G$ as in equation \eqref{Snegra}, then the
natural transformations $Nat(\bd{\ring{S}},\bd{\ring{S}'})$ form a
set. Moreover, there is a bijection
$$\hom{{}_{\bd{F}}\ring{C}_{\bd{G}}}{\,(S,\msf{s})}{\,(S',\msf{s}')\,}
\,\, \simeq \,\, Nat(\bd{\ring{S}},\bd{\ring{S}'}),$$
explicitly given by $$ \left( \xymatrix{ [\alpha: SF \rightarrow S']
\ar@{|->}[r] & [\bd{\alpha}: \bd{\ring{S}} \rightarrow
\bd{\ring{S}'}] } \right),\quad \left( \xymatrix{ S'\xi_{-}
\,\circ\, \bd{\beta}_{T^F(-)} & \bd{\beta} \ar@{|->}[l] }\right),$$
where for every object $(X,d^X) \in \cat{A}^F$, $\bd{\alpha}_X =
\alpha_X \,\circ\, Sd^X: S(X) \rightarrow S'(X)$.
\end{proposition}
\begin{proof}
We first prove that $Nat(\bd{\ring{S}},\bd{\ring{S}'})$ is a set. To
do this, we follow the proof of \cite[Lemma
5.1]{Iglesias/Gomez/Nastasescu:1999}. Let $\bd{\alpha},\bd{\beta}:
\bd{\ring{S}} \rightarrow \bd{\ring{S}'}$ two natural
transformations and $U$ a generator of $\cat{A}$. So
$(F(U),\delta_U)$ is a sub-generator of $\cat{A}^F$, by Lemma
\ref{additiva}. We claim that if $\bd{\alpha}_{(F(U),\,\delta_U)}
=\bd{\beta}_{(F(U),\,\delta_U)}$ then $\bd{\alpha}=\bd{\beta}$.
Considering any object $(X,d^X)$ of $\cat{A}^F$ with an epimorphism $\pi:
U^{(I)} \rightarrow X \rightarrow 0$ in $\cat{A}$, for some set $I$,
we obtain a diagram
$$ \xymatrix@C=40pt{ F(U^{(I)}) \cong F(U)^{(I)} \ar@{->}^-{\pi'}[r] &
F(X) \ar@{->}[r] & 0 \\ & X \ar@{->}_-{d^X}[u] & }$$ of morphisms of
$\cat{A}^F$. Since $\bd{\alpha}_{\left(F(U)^{(I)},\,d^{F(U)^{(I)}}
\right)} =\bd{\beta}_{\left(F(U)^{(I)},\,d^{F(U)^{(I)}}\right)}$, we
have
$$\bd{\ring{S}'}(\pi')\left(\underset{}{}  \bd{\alpha}_{\left(F(U)^{(I)},\,d^{F(U)^{(I)}} \right)}
- \bd{\beta}_{\left(F(U)^{(I)},\,d^{F(U)^{(I)}}\right)} \right) \,\,
=\,\, \left(\underset{}{}  \bd{\alpha}_{(F(X),\,\delta_X)} -
\bd{\beta}_{(F(X),\,\delta_X)}  \right) \bd{\ring{S}}(\pi')\,\,
=\,\,0.$$ By hypothesis, Lemma \ref{additiva} and diagram
\eqref{Snegra1}, we know that $\bd{\ring{S}}$ preserves
epimorphisms. Therefore, $\bd{\alpha}_{(F(X),\,\delta_X)} =
\bd{\beta}_{(F(X),\,\delta_X)}$, as $\pi'$ is an epimorphism. This
implies that $$\bd{\ring{S}'}(d^X) \,\,\circ\, \,
\left(\bd{\alpha}_{(X,\,d^X)} - \bd{\beta}_{(X,\,d^X)} \right)\,\,
=\,\, \left( \bd{\alpha}_{(F(X),\,\delta_X)} -
\bd{\beta}_{(F(X),\,\delta_X)} \right) \,\,\circ\, \,
\bd{\ring{S}}(d^X) \,\, =\,\, 0.$$ Applying the functor $S'^G$ and
using the diagram \eqref{Snegra1} for $(S',\sf{s}')$, we get $$ S'
S^F(d^X) \,\,\circ\, S'^G\left(\bd{\alpha}_{(X,\,d^X)} -
\bd{\beta}_{(X,\,d^X)} \right)\,\, =\,\, 0.$$ Composing with the map
$S'(\xi_X)$, we obtain that $S'^G\left(\bd{\alpha}_{(X,\,d^X)} -
\bd{\beta}_{(X,\,d^X)} \right) =0$ in $\cat{B}$, that is,
$\bd{\alpha}_{(X,\,d^X)} = \bd{\beta}_{(X,\,d^X)}$ in $\cat{B}^G$,
and this proves the claim.

For the stated bijection, we only prove that the mutually inverse
maps are well defined. Starting with a $2$-cell $\alpha: SF \rightarrow
S'$, and taking an arbitrary object $(X,d^X) \in \cat{A}^F$, we need
to demonstrate that $\bd{\alpha}_X =\alpha_X \,\circ\, Sd^X$ is a morphism
of the category $\cat{B}^G$, so
\begin{eqnarray*}
  d^{S'(X)} \,\circ\, \alpha_X \,\circ\, Sd^X &=& \msf{s}'_X \,\circ\, S'd^X \,\circ\, \alpha_X \,\circ\, Sd^X \\
   &=& \msf{s}'_X \,\circ\, \alpha_{F(X)} \,\circ\, SFd^X \,\circ\, Sd^X \\
   &=& \LR{\msf{s}'_X \,\circ\, \alpha_{F(X)} \,\circ\, S\delta_X} \,\circ\, Sd^X, \quad \text{ by }\eqref{coacion} \\
   &=& G\alpha_X \,\circ\, \msf{s}_{F(X)} \,\circ\, S\delta_X \,\circ\, Sd^X, \quad \text{ by }\eqref{dos-cells} \\
   &=& G\alpha_X \,\circ\, \msf{s}_{F(X)} \,\circ\, SFd^X \,\circ\, Sd^X \\
   &=& G\alpha_X \,\circ\, GSd^X \,\circ\, \msf{s}_X \,\circ\, Sd^X, \quad \msf{s}_{-} \text{ is natural} \\
   &=& G(\alpha_X \,\circ\, Sd^X) \,\circ\, d^{S(X)}.
\end{eqnarray*}
Obviously $\bd{\alpha}_{-}$ is natural. Conversely, starting with a
natural transformation $\bd{\beta}: \bd{\ring{S}} \rightarrow
\bd{\ring{S}'}$, its image is the natural transformation $\beta: SF
\rightarrow S'$ defined in every object $Y \in \cat{A}$ by
$\beta_{Y}= S'\xi_Y \,\circ\, \bd{\beta}_{T^F(Y)}$; we need to
show that $\beta$ satisfies the $2$-cell condition. In one hand, we have
\begin{eqnarray*}
  G \beta_{Y}  \,\circ\, \msf{s}_{F(Y)} \,\circ\, S\delta_{Y} &=& G S'\xi_Y \,\circ\, G\bd{\beta}_{T^F(Y)}
  \,\circ\, \msf{s}_{F(Y)} \,\circ\, S\delta_{Y}\\
   &=& G S'\xi_Y \,\circ\, G\bd{\beta}_{T^F(Y)} \,\circ\, d^{SF(Y)} \\
   &=& G S'\xi_Y \,\circ\, d^{S'F(Y)} \,\circ\, \bd{\beta}_{T^F(Y)}, \quad \bd{\beta}_{T^F(Y)}
   \text{ satisfies }\eqref{colineal} \\
   &=& G S'\xi_Y \,\circ\, \msf{s}'_{F(Y)} \,\circ\, S'\delta_Y \,\circ\, \bd{\beta}_{T^F(Y)} \\
   &=& G S'\xi_Y \,\circ\, \msf{s}'_{F(Y)} \,\circ\, \bd{\beta}_{T^F(F(Y))} \,\circ\, S\delta_Y,
   \quad  \bd{\beta}_{-} \text{ is natural} \\
   &=& \msf{s}'_Y \,\circ\, S'F\xi_Y \,\circ\, \bd{\beta}_{T^F(F(Y))} \,\circ\, S\delta_Y,
   \quad  \msf{s'}_{-} \text{ is natural} \\
   &=& \msf{s}'_Y \,\circ\, \bd{\beta}_{T^F(Y)} \,\circ\, SF\xi_Y \,\circ\,
   S\delta_Y ,
   \quad  \bd{\beta}_{-} \text{ is natural} \\ &=& \msf{s}'_Y \,\circ\, \bd{\beta}_{T^F(Y)}.
\end{eqnarray*}
On the other hand, we have
\begin{eqnarray*}
  \msf{s}'_Y \,\circ\, \beta_{F(Y)} \,\circ\, S\delta_Y &=& \msf{s}'_Y \,\circ\, S'\xi_{F(Y)}
  \,\circ\, \bd{\beta}_{T^F(F(Y))} \,\circ\, S\delta_Y \\
   &=& \msf{s}'_Y \,\circ\, S'\xi_{F(Y)}
  \,\circ\, S'\delta_Y  \,\circ\, \bd{\beta}_{T^F(Y)} \\
   &=& \msf{s}'_Y \,\circ\,  \bd{\beta}_{T^F(Y)}.
\end{eqnarray*}Therefore, $G \beta_{Y}  \,\circ\, \msf{s}_{F(Y)} \,\circ\, S\delta_{Y} = \msf{s}'_Y \,\circ\, \beta_{F(Y)} \,\circ\,
S\delta_Y$, for every object $Y \in \cat{A}$, and this gives the
needed condition.
\end{proof}

\section{Comonads and corings over rings with local
units}\label{sectionI}

We will consider rings without identity although we assume that a set
of identities is given. Following  \cite{Abrams:1983} (see also
\cite{Anh/Marki:1983} and \cite{Anh/Marki:1987}) a $\Bk$--module $A$
is said to be a \emph{ring with local units} if for every $a_1,\cdots,a_n$
in $A$, there exists an idempotent element $e \in \idemp{A}$ (the set
of all idempotent elements) such that
\begin{eqnarray*}
  a_ie &=& ea_i \,\, = \,\, a_i, \quad i=1,\cdots,n.
\end{eqnarray*}
We say that $e$ is a \emph{unity} for the set $\{a_1,\cdots,a_n\}$.
This is equivalent to say that for every $a,a' \in A$, there exists
a ring with identity of the form $A_e=eAe$ for some idempotent
element $e \in A$ such that $a, a' \in A_e$. For instance, every
induced ring from a $\Bk$--additive small category is a ring with
local units, in such case it is a ring with enough orthogonal idempotents,
see \cite{Gabriel:1962} and \cite{Fuller/Hullinger:1978}.

For any right $A$--module $X$ (i.e., a $\Bk$-module $X$ with
associative $\Bk$-linear right $A$-action $\mu_X: X\tensor{\Bk}A
\rightarrow X$), $XA$ denotes the right $A$-submodule
$$ XA \,=\, \left\{ \sum_{1\leq i \leq n} x_ia_i|\,\, x_i \in X,\,
a_i \in A, \text{ and } n \in \naturales \right\}.$$

A morphism between two rings with local units is a morphism of rings
$\psi: B \rightarrow A$ (i.e., compatible with multiplications)
satisfying the following condition: For every $a \in A$, there
exists $f \in \idemp{B}$ such that $a\psi(f) = \psi(f)a=a$. Observe
that this condition is equivalent to say that for every $e \in
\idemp{A}$, there exists $f \in \idemp{B}$ such that $e \psi(f) =
\psi(f) e = e$.

The construction of the usual tensor product over rings with
identity can be directly transferred to rings with local units, and
the most useful properties of this product are preserved. We use the
same symbol $-\tensor{A}-$ to denote the tensor product between
$A$--modules and $A$--linear morphisms for any ring with local units
$A$.

Let $A$ be a ring with local units, and $e \in \idemp{A}$. The
underlying  $\Bk$--module of the right $A$--module $eA$ is a direct
summand of $A$ with decomposition $A=eA \oplus \lr{ a-ea|\, a \in
A}$. Associated to $eA$ there are two $\Bk$--linear natural
transformations
\begin{equation}\label{gama-tau}
\xymatrix@R=0pt{\gamma_{e,X}: X \ar@{->}[r] & eA \tensor{A}X, \\ x
\ar@{|->}[r] & e\tensor{A}x } \quad \xymatrix@R=0pt{ \tau_{e,X}: eA
\tensor{A}X \ar@{->}[r] & X \\ ea\tensor{A}x \ar@{|->}[r] & eax, }
\end{equation}
for every right $A$--module $X$. Moreover, if $X$ is an
$(A,B)$--bimodule ($B$ is another ring with local units), then
$\gamma_{e,X}$ and $\tau_{e,X}$ are clearly right $B$--linear.
Taking $e' \in \idemp{A}$ another idempotent and $f: eA \rightarrow
e'A$ a right $A$--linear map, there are two commutative diagrams
\begin{equation}\label{gama-tau1}
\xymatrix@R=30pt@C=50pt{ X \ar@{->}^-{\gamma_{e,X}}[r]
\ar@{->}_-{\lambda_{f(e)}}[d] & eA \tensor{A}X
\ar@{->}^-{f\tensor{A}X}[d]
\\ X \ar@{->}^-{\gamma_{e',X}}[r] & e'A\tensor{A}X, } \qquad \xymatrix@R=30pt@C=50pt{
eA\tensor{A}X \ar@{->}^-{\tau_{e,X}}[r] \ar@{->}_-{f\tensor{A}X}[d]
& X \ar@{->}^-{\lambda_{f(e)}}[d] \\
e'A\tensor{A}X \ar@{->}^-{\tau_{e',X}}[r] & X }
\end{equation}
where $\lambda_{f(e)}: X \rightarrow X$ is the left multiplication
by $f(e)$, and $X$ any right $A$--module. Following
\cite{Abrams:1983}, there is a partial ordering on $\idemp{A}$ defined by
$$ e \leq e' \Longleftrightarrow e=ee'=e'e $$ for every $e,e' \in
\idemp{A}$. Taking $X_A$ any right $A$--module, and $e, e' \in
\idemp{A}$ such that $e \leq e'$, we can define a canonical
injection $\mu_{ee'}: Xe \rightarrow Xe'$, $\mu_e: Xe \rightarrow
X$. Furthermore, if $e \leq e' \leq e''$, then it is clear that
$\mu_{ee''} = \mu_{e'e''} \,\circ\, \mu_{ee'}$, thus
$\{(Xe,\mu_e)\}_{e \in \idemp{A}}$ is a directed system of
$\Bk$--submodule of $X$. In this way it is obvious that ${}_AA =
\dlimit(Ae)$ and $A_A= \dlimit(eA)$

A right $A$--module $X$ is said to be \emph{unital} if $XA=X$ (or $X
\cong X \tensor{A} A$ as right $A$--module, where the isomorphism
should be given by the right $A$--action). Equivalently, for every
element $m \in M $, there exists $e \in \idemp{A}$ such that $me=m$.
We denote by $\rmod{A}$ the full subcategory of the category of
right $A$--modules whose objects are all unital right $A$--modules.
An easy argument shows that $XA$ is the largest unital right
$A$--submodule of the right $A$--module $X$. On the other hand, a
right $A$--module $X$ is unital if and only if $\dlimit(Xe) =X$ in
the category of $\Bk$--modules. Given $B$ another ring with local
units, an \emph{unital} $(B,A)$--\emph{bimodule} is a
$(B,A)$--bimodule which is unital as a left $B$--module and as a
right $A$--module. Over the same ring, an $A$--bimodule $X$ is
unital if and only if for every $x \in X$, there exists $e
\in\idemp{A}$ such that $ex=xe=x$. In this way a morphism of rings
with local units $\psi: B \rightarrow A$ induces a structure of an
unital $B$--bimodule over $A$, and preserves the usual adjunction
between the categories of unital right modules
$$\xymatrix@C=60pt{ -\tensor{B}A: \rmod{B} \ar@<0,5ex>[r] &
\ar@<0,5ex>[l] \rmod{A}: \mathcal{O}. }$$

The definition of corings over rings with identity as it was introduced
in \cite{Sweedler:1975} can be directly extended, using unital
bimodules, to rings with local units. Let $A$ be a ring with local
units, an $A$--\emph{coring} is a three-tuple $(\cC,\Delta_{\cC},
\varepsilon_{\cC})$ consisting of an unital $A$--bimodule $\cC$ and
two $A$--bilinear maps
$$\xymatrix@C=50pt{\coring{C} \ar@{->}^-{\Delta_{\coring{C}}}[r] & \coring{C}\tensor{A}\coring{C}},\quad
\xymatrix@C=30pt{ \coring{C} \ar@{->}^-{\varepsilon_{\coring{C}}}[r]
& A}$$ such that $(\Delta_{\coring{C}}\tensor{A}\coring{C})
\,\circ\, \Delta_{\coring{C}} =
(\coring{C}\tensor{A}\Delta_{\coring{C}}) \,\circ\,
\Delta_{\coring{C}}$ and
$(\varepsilon_{\coring{C}}\tensor{A}\coring{C}) \,\circ\,
\Delta_{\coring{C}}=(\coring{C}\tensor{A}\varepsilon_{\coring{C}})
\,\circ\, \Delta_{\coring{C}}= \coring{C}$. A morphism of
$A$--corings is an $A$--bilinear map $\phi: \coring{C} \rightarrow
\coring{C}'$ which satisfies $\varepsilon_{\coring{C}'} \,\circ\,
\phi = \varepsilon_{\coring{C}}$ and $ \Delta_{\coring{C}'}
\,\circ\, \phi = (\phi \tensor{A} \phi) \,\circ\,
\Delta_{\coring{C}}$. We denote by $\coanillos{A}$ the category of
all $A$--corings and their morphisms.

A right $\coring{C}$--comodule is a pair $(M,\rho_{M})$ consisting
of an unital right $A$--module $M$ and a right $A$--linear map
$\rho_{M}: M \rightarrow M\tensor{A}\coring{C}$, called right
$\coring{C}$--coaction, such that $(M\tensor{A}\Delta_{\coring{C}})
\,\circ\, \rho_M = (\rho_M\tensor{A}\coring{C}) \,\circ\, \rho_M$
and $(M\tensor{A}\varepsilon_{\coring{C}}) \,\circ\, \rho_M=M$. A
morphism of right $\coring{C}$--comodules is a right $A$--linear map
$f: M \rightarrow M'$ satisfying $ \rho_{M'} \,\circ\, f =
(f\tensor{A}\coring{C}) \,\circ\, \rho_M$. Right
$\coring{C}$--comodules and their morphisms form a not necessarily
abelian category which we denote by $\rcomod{\coring{C}}$ (see
\cite[Section 1]{Kaoutit/Gomez/Lobillo:2004c}). For every unital
right $A$--module $X$ the pair
$(X\tensor{A}\cC,X\tensor{A}\Delta_{\cC})$ is clearly a right
$\cC$--comodule. This establishes in fact a functor $-\tensor{A}\cC:
\rmod{A} \rightarrow \rcomod{\cC}$ with the forgetful functor
$U_{\cC}: \rcomod{\cC} \rightarrow \rmod{A}$ as a left adjoint (see
\cite{Guzman:1989}).

\begin{example}\label{corings}
Of course every ring with local units $A$ is trivially an
$A$--coring with comultiplication the isomorphism $A \cong
A\tensor{A} A$ and counit the identity map $A$.
\begin{enumerate}
\item \cite{Sweedler:1975}. Let $\psi: B \rightarrow A$ be a
morphism of rings with local units and consider the unital
$A$--bimodule $A\tensor{B}A$ with the following two maps
$$\xymatrix@R=0pt{\Delta: A \tensor{B} A \ar@{->}[r] & (A \tensor{B}A)
\tensor{A} (A \tensor{B} A), \\ a \tensor{B} a' \ar@{|->}[r] & a
\tensor{B}e \tensor{A}e \tensor{B}a',} \quad
\xymatrix@R=0pt{\varepsilon : A\tensor{B}A \ar@{->}[r] & A \\ a
\tensor{B} a' \ar@{|->}[r] & aa', }
$$ where $e \in \psi(\idemp{B})$ such that $ea=ae=a$ and
$a'e=ea'=a'$; that is, $e$ is a unity for both $a$ and $a'$. To check
that $\Delta(a\tensor{B}a')$ is independent of the choice of the
unity, let us consider another unity $e' \in \psi(\idemp{B})$ for
both $a$ and $a'$. By definition there exists $e'' \in
\psi(\idemp{B})$ a unity for $e$ and $e'$. Of course $e''$ is also a
unity for both $a$ and $a'$. Now, we have
\begin{eqnarray*}
  a\tensor{B}e\tensor{A}e \tensor{B}a' &=& a\tensor{B}ee''\tensor{A}e''e \tensor{B}a',
  \quad \text{ since }\, e=ee''=e''e  \\
   &=& ae\tensor{B}e''\tensor{A}e'' \tensor{B}ea' \\
   &=& a\tensor{B}e''\tensor{A}e'' \tensor{B}a'.
\end{eqnarray*}
Similarly, we get $a\tensor{B}e'\tensor{A}e' \tensor{B}a'\,=\,
a\tensor{B}e''\tensor{A}e'' \tensor{B}a'$, and so
$\Delta(a\tensor{B}a')$ is independent from the choice of $e$. An
easy verification shows now that $(A\tensor{B}A,
\Delta,\varepsilon)$ is an $A$--coring.

\item Let $M$ be an unital
$A$--bimodule over a ring with local units $A$. Consider the direct
sum of an $A$--bimodules $\cC: = A \oplus M$ together with the
$A$--bilinear maps
$$\xymatrix@R=0pt{\Delta: \coring{C} \ar@{->}[r] &
\coring{C}\tensor{A}\coring{C} \\ (a,m) \ar@{|->}[r] &
(a,0)\tensor{A}(e,0) + (0,m)\tensor{A}(e,0) \\ & +
(e,0)\tensor{A}(0,m) }\quad \xymatrix@R=0pt{ \varepsilon: \coring{C}
\ar@{->}[r] & A \\ (a,m) \ar@{|->}[r] & a }
$$ where $e \in \idemp{A}$ such that $em=me=m$ and $ea=ae=a$. Let
us check that $\Delta$ is a well defined map. First we observe that a
common unity for $a$ and $m$ does always exist. If $e' \in \idemp{A}$
is another unity for $a$ and $m$, then one can consider $e'' \in
\idemp{A}$ a unity for $e'$ and $e$. Therefore, we have
\begin{eqnarray*}
 \Delta(a,m) &=& (a,0)\tensor{A}(e,0) + (0,m)\tensor{A}(e,0) +
(e,0)\tensor{A}(0,m) \\ &=& (a,0)\tensor{A}(ee'',0) +
(0,m)\tensor{A}(ee'',0) +
(e''e,0)\tensor{A}(0,m) \\
   &=& (a,0)e\tensor{A}(e'',0) + (0,m)e\tensor{A}(e'',0) +
(e'',0)\tensor{A}e(0,m) \\
   &=& (a,0)\tensor{A}(e'',0) + (0,m)\tensor{A}(e'',0) +
(e'',0)\tensor{A}(0,m).
\end{eqnarray*}
In the same way, we get $(a,0)\tensor{A}(e'',0) +
(0,m)\tensor{A}(e'',0) + (e'',0)\tensor{A}(0,m)\,=\,
(a,0)\tensor{A}(e',0) + (0,m)\tensor{A}(e',0) +
(e',0)\tensor{A}(0,m)$, and thus $\Delta(a,m)$ is independent from
the choice of the unity $e$. The three-tuple
$(\cC,\Delta,\varepsilon)$ is easily proved to be an $A$--coring.

\item \cite{Kaoutit/Gomez:2003a}. Let ${}_B\Sigma_A$ be an unital bimodule
over rings with local units $A$ and $B$ such that $\Sigma_A$ is a
finitely generated and projective unital right module with finite right
dual basis $\{(u_i,u_i^*)\}_i \subset \Sigma \times
\Sigma^*$ where $\Sigma^*=\hom{A}{\Sigma}{A}$. That is, $u \,=\,
\sum_iu_iu_i^*(u)$, for every $u \in \Sigma$.  It is well known that
$\Sigma^*$ is also an unital $(A,B)$--bimodule, and thus $\Sigma^*
\tensor{B} \Sigma$ is an unital $A$--bimodule. Furthermore, there
exist two $A$--bilinear maps
$$\xymatrix@R=0pt{ \Delta: \Sigma^* \tensor{B} \Sigma \ar@{->}[r]
& \Sigma^* \tensor{B} \Sigma \tensor{A} \Sigma^* \tensor{B} \Sigma,
\\ u^* \tensor{B} u \ar@{|->}[r] & \sum_{i} u^* \tensor{B} u_i
\tensor{A} u^*_i \tensor{B} u } \quad \xymatrix@R=0pt{ \varepsilon:
\Sigma^* \tensor{B} \Sigma \ar@{->}[r] & A \\ u^* \tensor{B} u
\ar@{|->}[r] & u^*(u) }$$ The canonical isomorphism $\Sigma
\tensor{A} \Sigma^* \cong \rend{\Sigma}{A}$ implies that $\Delta$ is
independent from the choice of this right dual basis, and that
$(\Sigma^* \tensor{B} \Sigma, \Delta,\varepsilon)$ is an
$A$--coring. This coring is known as \emph{the finite comatrix
coring} associated to ${}_{B}\Sigma_A$.
\end{enumerate}
\end{example}

\begin{example}
To take a specific example of finite comatrix corings over rings
with local units, we consider the so called \emph{finitely
orthogonal Rees matrix rings} extensively investigated in
\cite{Anh/Marki:1983}. Following \cite[Example 2]{Anh/Marki:1987}
(see also \cite{Anh/Marki:1983} for notions occurring here), let $R$
be a ring with identity, and $A$ a Rees matrix ring over $R$ with
canonical decompositions $A \cong Ae \tensor{eAe} eA$, $e \in
\idemp{A}$ and $R \cong eAe$. If $A$ is left-right finitely
orthogonal with respect to $e$ \cite[Definition
4.2]{Anh/Marki:1983}, then one can easily prove that $A$ is a ring
with local units. On the other hand, if we take
${}_{eAe}\Sigma_A=eA$, then the associated finite comatrix
$A$--coring is given by the $A$--bimodule $Ae\tensor{eAe}eA$, and
its counit is just the above isomorphism $\xymatrix{Ae\tensor{eAe}eA
\ar@{->}^-{\cong}[r] &  A}$ sending $ae\tensor{eAe}ea' \mapsto
aea'$. Therefore, we have an isomorphism of categories
$\rcomod{Ae\tensor{eAe}eA} \cong \rmod{A}$ via this counit (recall
that the counit is always a morphism of corings). Since the right
$(Ae\tensor{eAe}eA)$--comodule $\Sigma$ is clearly a generator of
$\rcomod{Ae\tensor{eAe}eA}$, we deduce following Gabriel-Popescu's
Theorem \cite{Gabriel/Popescu:1964} that ${}_{eAe}\Sigma$ is a
faithfully flat module (here $eAe$ coincides with the endomorphism
ring of this comodule). Thus, $-\tensor{eAe}\Sigma_A: \rmod{eAe}
\rightarrow \rcomod{Ae\tensor{eAe}eA}$ establishes an equivalence of
categories by using the non unital version of the generalized
Descent Theorem \cite[Theorem 3.10]{Kaoutit/Gomez:2003a}. In
conclusion, $A$ is Morita equivalent to $eAe$, and thus to $R$,
which gives an alternative proof of \cite[Example
2]{Anh/Marki:1987}.
\end{example}

From now on, we fix $A$, $B$ rings with local units. Let $F :
\rmod{A} \rightarrow \rmod{B}$ be a continuous functor
(thus it also preserves inductive limits). As in the case of rings with identity
\cite{Watts:1960}, next we will show that $F$ is naturally
isomorphic to the tensor product functor. Another approach,
concerning functors valued in abelian groups was given in
\cite{Fisher/Newell:1971}.

The structure of an $(A,B)$--bimodule over $F(A)$ comes out by the
composition map
\begin{equation}\label{biaction}
\xymatrix@C=50pt{ A \ar@{->}^-{\lambda}[r] & \hom{A}{A_A}{A_A}
\ar@{->}^-{F}[r] & \hom{B}{F(A)}{F(A)} }
\end{equation}
where $\lambda_a: A_A \rightarrow A_A$ is the left multiplication by
$a \in A$. In the same way we get an $(A,B)$--biaction over
$F^n(A) = (F \,\circ\, \cdots \,\circ\, F)(A)$ ($n$-times).
Therefore, one can consider the functor $-\tensor{A}F(A): \rmod{A}
\rightarrow \rmod{B}$. Now, let's start with an arbitrary idempotent
element $e \in \idemp{A}$, and consider the composed right
$B$--linear map
\begin{equation}\label{Ups}
\xymatrix@R=20pt@C=60pt{ F(eA) \ar@{-->}^-{\Upsilon_{eA}}[rr]
\ar@{->}_-{F(\tau_{e})}[dr]
& & eA\tensor{A}F(A) \\
& F(A) \ar@{->}_-{\gamma_{e,F(A)}}[ur] &  }
\end{equation}
where $\tau_e: eA \rightarrow A$ is the canonical injection, and
$\gamma_{e,F(A)}$ is the map defined in equation
\eqref{gama-tau}. If we take $e' \in \idemp{A}$ another idempotent
element and $f: eA \rightarrow e'A$ a right $A$--linear map, then by
equation \eqref{gama-tau1}, we get a commutative diagram like this
$$\xymatrix{F(eA) \ar@{-->}^-{\Upsilon_{eA}}[rr]
\ar@{->}_-{F(\tau_{e})}[dr] \ar@{->}_-{F(f)}[dd]
& & eA\tensor{A}F(A) \ar@{->}^-{f\tensor{A}F(A)}[dd] \\
& F(A) \ar@{->}_-{\gamma_{e,F(A)}}[ur] \ar@{->}^<<{F(\lambda_{f(e)})}[dd] & \\
F(e'A) \ar@{-->}^-{\Upsilon_{e'A}}[rr] \ar@{->}_-{F(\tau_{e'})}[dr]
& & e'A\tensor{A}F(A) \\
& F(A) \ar@{->}_-{\gamma_{e',F(A)}}[ur] &  }
$$
which shows that $\Upsilon_{-}$ is natural over the set of right
$A$--modules $\{eA\}_{e \in \idemp{A}}$. Using the projections
$\pi_e: A \rightarrow eA$ and the maps $\tau_{e,F(A)}$ defined in
equation \eqref{gama-tau}, we can also construct a right $B$--linear
map
\begin{equation}\label{Theta}
\xymatrix@R=20pt@C=60pt{ eA\tensor{A}F(A)
\ar@{-->}^-{\Theta_{eA}}[rr] \ar@{->}_-{\tau_{e,F(A)}}[dr]
& & F(eA) \\
& F(A) \ar@{->}_-{F(\pi_e)}[ur] &  }
\end{equation}
which is natural over $\{eA\}_{e \in \idemp{A}}$ by equation
\eqref{gama-tau1}.

\begin{lemma}\label{ISO}
Let $F: \rmod{A} \rightarrow \rmod{B}$ be a continuous functor. For every idempotent
$e \in \idemp{A}$, we have
\begin{eqnarray*}
  \Upsilon_{eA} \, \,\circ\, \, \Theta_{eA} &=& eA \tensor{A} F(A),  \\
  \Theta_{eA} \, \,\circ\, \, \Upsilon_{eA}&=& F(eA).
\end{eqnarray*}
In particular $\Upsilon_{-}$ extended uniquely to a natural
isomorphism $$\xymatrix{ \Upsilon_{-}^F : F(-) \ar@{->}[r] &
-\tensor{A}F(A), }$$ and $F(A)$ becomes a left unital $A$--module,
and thus an unital $(A,B)$--bimodule.
\end{lemma}
\begin{proof}
By definition, the left $A$--action of $F(A)$ is given by the rule
\begin{eqnarray*}
  a \, x &=& F(\lambda_a)(x), \, \text{ for every pair } (a,x)
  \in A \times F(A).
\end{eqnarray*}
Fixing an arbitrary idempotent $e \in\idemp{A}$. So \begin{eqnarray*}
  \Upsilon_{eA} \,\circ\, \Theta_{eA} (ea\tensor{A}x) &=& \Upsilon_{eA} \,\circ\, F(\pi_e) (ea \, x) \\
   &=& \Upsilon_{eA} \,\circ\, F(\pi_e) \,\circ\, F(\lambda_{ea}) ( x) \\
   &=& \Upsilon_{eA} \,\circ\, F(\lambda_{ea}) (x) \\
   &=& \gamma_{e,F(A)} \,\circ\, F(\tau_e \,\circ\, \lambda_{ea}) (x) \\
   &=& \gamma_{e,F(A)} \,\circ\, F(\lambda_{ea}) (x) \\
   &=& \gamma_{e,F(A)} (ea \, x) \,\, = \,\, e\tensor{A} eax, \\
   &=& ea\tensor{A}x
\end{eqnarray*}
for every $a \in A$, $x \in F(A)$, and this shows the first
equality. To check the second equality, take $y \in F(eA)$ and
compute
\begin{eqnarray*}
  \Theta_{eA} \,\circ\, \Upsilon_{eA}(y) &=& \Theta_{eA} \,\circ\, \gamma_{e,F(A)} \,\circ\, F(\tau_e) (y) \\
   &=& F(\pi_e) \,\circ\, \tau_{e,F(A)} \LR{ e \tensor{A} F(\tau_e)(y) } \\
   &=& F(\pi_e) \,\circ\, \LR{ e \, F(\tau_e)(y) } \\
   &=& F(\pi_e) \,\circ\, \LR{ F(\lambda_e) \,\circ\, F(\tau_e) (y)} \\
   &=& F( \pi_e) \,\circ\, \LR{F(\lambda_e \,\circ\, \tau_e) (y)} \\
   &=& F(\pi_e) \,\circ\, F(\tau_e) (y) \,\, = \,\, y.
\end{eqnarray*}
$\Upsilon_{-}$ is clearly extended to unital right $A$--modules of
the form $\oplus_{j \in J} (e_jA)^{(I_j)}$ ($J$ and $I_j$ are sets).
Since $F$ preserves direct sums, this extension is also an
isomorphism of unital right $A$--modules. By Mitchell's
Theorem \cite[Theorem 4.5.2]{Mitchell:1965} (\cite[Theorem
3.6.5]{Popescu:1973}), $\Upsilon_{-}$ extends uniquely to a natural
isomorphism $\Upsilon_{-}^F: F(-) \rightarrow -\tensor{A}F(A)$ over
all unital right $A$--modules.

We need to show that $F(A)$ is a left unital $A$--module. We have
seen that $eF(A) \cong F(eA)$ for every idempotent $e \in \idemp{A}$
via $\Upsilon_{eA}^F$. Since $F$ preserves inductive limits, we get
$$\dlimit{(eF(A))} \,\cong \, \dlimit{(F(eA))} \,\cong \,F(\dlimit{(eA)}) \,\cong \,F(A)$$
which implies that ${}_AF(A)$ is left unital and this finishes the
proof.
\end{proof}

The following lemma can be deduced from
\cite[39.5]{Brzezinski/Wisbauer:2003}. For the sake of completeness,
we include a detailed proof.

\begin{lemma}\label{comonad1}
Let $A$, $B$ and $C$ be three rings with local units. We denote by
$\Upsilon_{-}^{\chi}$ the natural isomorphism of Lemma \ref{ISO}
associated to the continuous functor $\chi$.
\begin{enumerate}[ ]
\item [(a)] Let $F_1,\,F_2: \rmod{A} \rightarrow \rmod{B}$ be two continuous functors. Assuming that there exists a
natural transformation $\Xi: F_1 \rightarrow F_2$, then
\begin{enumerate}
\item[(1)] The morphism $\Xi_A: F_1 (A) \rightarrow F_2(A)$ is
$(A,B)$--bilinear. \item[(2)] For every unital right $A$--module
$X$, we have a commutative diagram $$\xymatrix@R=30pt@C=60pt{ F_1(X)
\ar@{->}^-{\Upsilon_X^{F_1}}[r]
\ar@{->}_{\Xi_X}[d] & X\tensor{A}F_1(A) \ar@{->}^-{X\tensor{A}\Xi_A}[d] \\
F_2(X) \ar@{->}^-{\Upsilon_X^{F_2}}[r] & X\tensor{A}F_2(A). }$$
\end{enumerate}
\item[(b)] Let $F:\rmod{A} \to \rmod{B}$ and $S:\rmod{B}\to
\rmod{C}$ be two continuous functors.
Then we have  $\Upsilon_{-}^{F(-)\tensor{B}S(B)} = \Upsilon_{-}^{F}
\tensor{B}S(B)$ and the following diagram of natural transformations
commutes $$\xymatrix@R=35pt@C=70pt{SF(-)
\ar@{->}|-{\Upsilon_{F(-)}^S}[d] \ar@{->}^-{\Upsilon_{-}^{SF}}[r]
&  -\tensor{A}SF(A) \ar@{->}|-{-\tensor{A}\Upsilon_{F(A)}^S}[d] \\
F(-)\tensor{B}S(B) \ar@{->}^-{\Upsilon_{-}^{F(-)\tensor{B}S(B)}}[r]
& -\tensor{A}F(A)\tensor{B}S(B) }$$
\end{enumerate}
\end{lemma}
\begin{proof} $(a)$ $(1)$. By definition, we only need to show that $\Xi_A$ is left
$A$--linear. For every element $a \in A$, we have $\Xi_A \,\circ\,
F_1(\lambda_a) = F_2(\lambda_a) \,\circ\, \Xi_A$ since $\Xi_{-}$ is
natural and $\lambda_a$ is right $A$--linear. Hence $\Xi_A (ax) = a
\Xi_A(x)$,
for every pair of elements $(a,x) \in A \times F_1(A)$.  \\
$(a)$ $(2)$. Assuming that $X$ is of the form $X=eA$ for some
idempotent element $e \in \idemp{A}$, we obtain
\begin{eqnarray*}
  \Upsilon_{eA}^{F_2} \, \,\circ\, \, \Xi_{eA} &=& \gamma_{e,F_2(A)} \,\circ\, F_2(\tau_e) \,\circ\, \Xi_{eA} \\
   &=& \gamma_{e,F_2(A)} \,\circ\, \Xi_A \,\circ\, F_1(\tau_e), \quad \Xi_{-} \text{ is natural} \\
   &=& \LR{eA\tensor{A}\Xi_A} \,\circ\, \gamma_{e,F_1(A)} \,\circ\,
   F_1(\tau_e), \quad
    \gamma_{e,-} \text{ is natural } \\
   &=& \LR{eA\tensor{A}\Xi_A} \,\circ\, \Upsilon_{eA}^{F_1}.
\end{eqnarray*}
For the general case we use a free presentation $\oplus_{k}
(e_{k}A)^{(I_{k})} \rightarrow X \rightarrow 0$ where $\{e_{k}\}_{k}
\subseteq \idemp{A}$, $I_{k}$ are sets, and the previous case taking
into account the hypothesis done over the stated functors. \\ $(b)$
Straightforward.
\end{proof}

It is clear that any $A$--coring $(\coring{C},\Delta,\varepsilon)$
induces, by the three-tuple
$((-\tensor{A}\coring{C}),-\tensor{A}\Delta,-\tensor{A}\varepsilon)$,
a comonad in $\rmod{A}$. Now, consider $(F,\delta,\xi)$ a comonad in
$\rmod{A}$ such that $F: \rmod{A} \rightarrow \rmod{A}$ is a continuous functor. Our next goal is to prove
that $F(A)$ admits a structure of an $A$--coring. Notice that this
has been already observed in \cite[p. 398]{Eilenberg/Moore:1965} (without
proofs) in the case of commutative rings with identity.

We begin by introducing some convenient notations. Denote by
$\Upsilon_{-}^{id_{\rmod{A}}}: id_{\rmod{A}} \rightarrow
-\tensor{A}A$ the canonical natural isomorphism, and by
$\Upsilon_{-}^{\chi}: \chi(-) \rightarrow -\tensor{A}\chi(A)$ the
natural isomorphism of Lemma \ref{ISO} associated to the continuous functor $\chi:\rmod{A} \rightarrow \rmod{A}$.
For every element $a \in A$,
$\lambda_a: A_A \rightarrow A_A$ stills denoting the left
multiplication by $a$.

\begin{proposition}\label{comonad2}
Let $A$ be a ring with local units, and $(F,\delta,\xi)$ a comonad
in $\rmod{A}$ such that $F$ is a continuous functor. Then
$(F(A),\Upsilon_{F(A)}^F \,\circ\, \delta_A, \xi_A)$ is an
$A$--coring.
\end{proposition}
\begin{proof}
By Lemma \ref{ISO}, $F(A)$ admits a structure of unital
$A$--bimodule. The maps $\Delta = \Upsilon_{F(A)}^F \,\circ\,
\delta_A$ and $\varepsilon = \xi_A$ are $A$--bilinear by
Lemma \ref{comonad1}(a)(1).\\
By hypothesis, we have the following diagram
$$\xymatrix@R=33pt@C=67pt{ \ar @{}[dr]|{(r1)} F(A)
\ar@{->}^-{\delta_A}[r] \ar@{->}|-{\delta_A}[d] & \ar @{}[dr]|{(r2)}
F^2(A) \ar@{->}^-{\Upsilon_{F(A)}^F}[r] \ar@{->}|-{F(\delta_A)}[d]
 & F(A)\tensor{A}F(A) \ar@{->}|-{\delta_A\tensor{A}F(A)}[d] \\
\ar @{}[dr]|{(r3)} F^2(A) \ar@{->}^-{\delta_{F(A)}}[r]
\ar@{->}|-{\Upsilon_{F(A)}^F}[d] & \ar @{}[dr]|{(r4)} F^3(A)
\ar@{->}^-{\Upsilon_{F^2(A)}^F}[r]
\ar@{->}|-{\Upsilon^{F^2}_{F(A)}}[d]
 & F^2(A)\tensor{A}F(A) \ar@{->}|-{\Upsilon_{F(A)}^F\tensor{A}F(A)}[d] \\
F(A)\tensor{A}F(A) \ar@{->}^-{F(A)\tensor{A}\delta_A}[r] & F(A)
\tensor{A}F^2(A) \ar@{->}^-{F(A)\tensor{A}\Upsilon_{F(A)}^F}[r] &
F(A)\tensor{A}F(A)\tensor{A}F(A) }$$ where the rectangle $(r1)$ is
commutative by equation \eqref{comonad}, and $(r2)$ by the
naturality of $\Upsilon_{-}^F$. Applying Lemma \ref{comonad1} to the
natural transformation $\delta: F \rightarrow F^2$, we get the
commutativity of the rectangle $(r3)$. Lastly, Lemma \ref{comonad1}
applied this time to $\Upsilon_{F(-)}^F: F^2(-) \rightarrow
F(-)\tensor{A}F(A)$, gives the commutativity of the rectangle
$(r4)$. Therefore, the total diagram is commutative; whence $\Delta$
is coassociative. The left counitary property is shown by the
following diagram
$$\xymatrix@R=30pt@C=60pt{ F(A) \ar@{->}^{\delta_A}[r] \ar@{=}[dr] & F^2(A)
\ar@{->}^-{\Upsilon_{F(A)}^F}[r] \ar@{->}|-{F(\xi_A)}[d]
& F(A)\tensor{A}F(A) \ar@{->}|-{\xi_A\tensor{A}F(A)}[d] \\
& F(A) \ar@{->}^-{\Upsilon_{A}^F}[r] & A\tensor{A}F(A) }$$ which is
commutative by the naturality of $\Upsilon_{-}^F$ and equation
\eqref{comonad}. Now, this last equation together with Lemma
\ref{comonad1} applied to $\xi: F \rightarrow id_{\rmod{A}}$ give
the commutative diagram
$$\xymatrix@R=30pt@C=60pt{ F(A) \ar@{->}^{\delta_A}[r] \ar@{=}[dr] & F^2(A)
\ar@{->}^-{\Upsilon_{F(A)}^F}[r] \ar@{->}|-{\xi_{F(A)}}[d]
& F(A)\tensor{A}F(A) \ar@{->}|-{F(A)\tensor{A}\xi_A}[d] \\
& F(A) \ar@{->}^-{\Upsilon_{F(A)}^{id_{\rmod{A}}}}[r] &
F(A)\tensor{A}A }$$ which leads to the right counitary property, and
this finishes the proof.
\end{proof}

Let $A$ be a ring with local units, and $(F,\delta,\xi)$ a comonad
in $\rmod{A}$. Consider the universal cogenerator of this comonad,
that is, an adjunction $\xymatrix{ S^F: \rmod{A}^F \ar@<0,5ex>[r] &
\ar@<0,5ex>[l] \rmod{A}: T^F}$ such that $F = S^F \,\circ\, T^F$,
and where $\rmod{A}^F$ is the category of comodules over
$(F,\delta,\xi)$. The functor $S^F$ is then the forgetful functor:
$S^F(X,d^X)=X$, for every comodule $(X,d^X)$, and the functor $T^F$
is defined over objects by $T^F(M) = (F(M), \delta_M)$, for every
right unital $A$--module $M$. Given an $A$--coring
$(\coring{C},\Delta,\varepsilon)$ it is easily seen that the
canonical adjunction $\xymatrix{ U_{\cC}: \rcomod{\coring{C}}
\ar@<0,5ex>[r] & \ar@<0,5ex>[l] \rmod{A}: -\tensor{A}\coring{C}}$ is
the universal cogenerator of the associated comonad $(
-\tensor{A}\coring{C}, -\tensor{A}\Delta, -\tensor{A}\varepsilon)$
in $\rmod{A}$. The following proposition compares the cogenerator of
the comonad $(F,\delta,\xi)$ with the canonical adjunction
associated to the coring $(F(A),\Upsilon_{F(A)}^A \circ
\delta_A,\xi_A)$ of Proposition \ref{comonad2}.

\begin{proposition}\label{comonad3}
Let $A$ be a ring with local units, and $(F,\delta,\xi)$ a comonad
in $\rmod{A}$  such that $F$ is a continuous functor with the
universal cogenerator $\xymatrix{ S^F: \rmod{A}^F \ar@<0,5ex>[r] &
\ar@<0,5ex>[l] \rmod{A}: T^F}$. Considering $F(A)$ as an $A$--coring
with the structure of Proposition \ref{comonad2}, we obtain an
isomorphism of categories
$$\xymatrix@C=60pt{ \digamma: \rmod{A}^F \ar@<0,5ex>^-{\cong}[r] &
\ar@<0,5ex>[l] \rcomod{F(A)} : \daleth,}$$ such that $ S^F \,\, =
\,\, U_{F(A)} \, \,\circ\, \, \digamma$,  and $ \digamma \,\circ\,
T^F \,\, \cong \,\, -\tensor{A}F(A)$ is a natural isomorphism.
\end{proposition}
\begin{proof}
Given $(X,d^X)$ an object of $\rmod{A}^F$, we get a diagram
\begin{equation}\label{diag1}
\xymatrix@R=30pt@C=70pt{ \ar@{}[dr]|{(r1)} X \ar@{->}^-{d^X}[r]
\ar@{->}|-{d^X}[d] & \ar@{}[dr]|{(r2)} F(X)
\ar@{->}^-{\Upsilon_X^F}[r] \ar@{->}|-{\delta_X}[d] & X \tensor{A}
F(A) \ar@{->}|-{X\tensor{A}\delta_A}[d] \\ \ar@{}[dr]|{(r3)} F(X)
\ar@{->}^-{F(d^X)}[r] \ar@{->}|-{\Upsilon_{X}^F}[d] &
\ar@{}[dr]|{(r4)} F^2(X) \ar@{->}^-{\Upsilon^{F^2}_X}[r]
\ar@{->}|-{\Upsilon_{F(X)}^F}[d] & X \tensor{A} F^2(X)
\ar@{->}|-{X\tensor{A}\Upsilon_{F(A)}^F}[d] \\
X \tensor{A}F(X) \ar@{->}^-{d^X\tensor{A}F(A)}[r] & F(X) \tensor{A}
F(A) \ar@{->}^-{\Upsilon_X^F\tensor{A}F(A)}[r] & X\tensor{A}F(A)
\tensor{A} F(A). }
\end{equation}
The rectangles $(r1)$ and $(r3)$ are by definition commutative.
Applying Lemma \ref{comonad1} consecutively  to the natural
transformations $\delta: F \rightarrow F^2$ and $\Upsilon_{F(-)}^F:
F^2(-) \rightarrow F(-)\tensor{A}F(A)$, we obtain the commutativity
of rectangles $(r2)$ and $(r4)$ (here the natural isomorphism
associated to the functor $F(-) \tensor{A}F(A)$ is
$\Upsilon_{-}^{F(-) \tensor{A}F(A)} =\Upsilon_{-}^F\tensor{A}F(A)$,
see Lemma \ref{comonad1}(b)). Thus, the whole diagram is commutative
which shows the coassociativity of the map $\Upsilon_X^F \,\circ\,
d^X$. We also have another commutative diagram
$$\xymatrix@=30pt@C=60pt{ X \ar@{->}^-{d^X}[r] & F(X)
\ar@{->}^-{\Upsilon_X^F}[r]
\ar@{->}|-{\xi_X}[d] & X\tensor{A}F(A) \ar@{->}|-{X\tensor{A}\xi_A}[d] & \\
& X \ar@{->}^-{\Upsilon^{id_{\rmod{A}}}_X}[r] & X\tensor{A}A
\ar@{->}^-{(\Upsilon^{id_{\rmod{A}}}_X)^{-1}}[r] & X. }$$ Therefore,
$(S^F(X)=X, \rho_X=\Upsilon_X^F \,\circ\, d^X)$ is a right
$F(A)$--comodule. Now any morphism $f: (X,d^X) \rightarrow
(X',d^{X'})$ in the category $\rmod{A}^F$ entails a commutative
diagram
$$\xymatrix@R=30pt@C=60pt{ X \ar@{->}^-{d^X}[r] \ar@{->}|-{f}[d] & F(X)
\ar@{->}^-{\Upsilon_X^F}[r] \ar@{->}|-{F(f)}[d] & X\tensor{A}F(A)
\ar@{->}|-{f\tensor{A}F(A)}[d] \\ X' \ar@{->}^-{d^{X'}}[r] & F(X')
\ar@{->}^-{\Upsilon_{X'}^F}[r]  & X'\tensor{A}F(A). }$$ Hence
$\digamma: \rmod{A}^F \rightarrow \rcomod{F(A)}$ defined by
$$\digamma(X,d^X) \, = \, (S^F(X), \Upsilon_X^F \,\circ\, d^X), \text{
and } \digamma(f)= S^F(f)$$ is a well defined functor with inverse
$\daleth: \rcomod{F(A)} \rightarrow \rmod{A}^F$ defined by
$$\daleth(Y,\rho_Y) \, = \, (Y, d^Y = (\Upsilon_Y^F)^{-1} \,\circ\,
 \rho_Y),  \text{ and } \daleth(g)=g.$$ Clearly the underlying
right $A$--module of $\digamma(X,d^X)$ coincides with that of
$(X,d^X)$. That is $U_{F(A)} \,\circ\, \digamma = S^F$. The
commutative rectangles $(r2)$ and $(r4)$ in diagram \eqref{diag1}
assert that $\Upsilon_X^F: \digamma \,\circ\, T^F(X) =
\digamma(F(X),\delta_X) \rightarrow \LR{X\tensor{A}F(A),
X\tensor{A}\left( \Upsilon_{F(A)}^F \,\circ\, \delta_A \right)}$ is
an isomorphism of right $F(A)$--comodules, for every right unital
module $X \in \rmod{A}$, and this leads to the stated natural
isomorphism $\Upsilon_{-}^F: \digamma \,\circ\, T^F(-) \cong
-\tensor{A}F(A)$.
\end{proof}

\begin{remark}\label{izq-dercha}
Of course one can work with left unital modules, and prove
similar results concerning functors which are right exact and
preserve direct sums and the induced corings from their comonads. In
this paper we only work with right unital modules. The left-right
relationship is omitted.
\end{remark}

\section{The bi-equivalence of bicategories}\label{section4}

In this section we define the bicategory of corings over rings with
local units in the same way as in
\cite{Brzezinski/Kaoutit/Gomez:2004}, using general methods from
\cite{Lack/Street:2002}. Next, we establish a bi-equivalence between
this bicategory and the $2$-category of comonads in right unital
modules over rings with local units as it was defined in Section
\ref{section0}.
\bigskip

In what follows $\mathscr{B}$ denotes the bicategory of unital bimodules
(i.e. $0$-cells are all rings  with local units,
Hom-Categories are categories of unital bimodules). The multiplications
are given by the tensor product bi-functors. Let us consider as in
Section \ref{section0}, the $2$-category $\td{\mathscr{L}}$ whose
$0$-cells are all Grothendieck categories of the form $\rmod{A}$ for
some ring with local units $A$, and whose Hom-Categories $\bara{
\rm Funct}(\rmod{A},\rmod{B})$ consists of continuous functors. The multiplications are given by the usual
compositions of functors and natural transformations. The identity
$1$-cell of a given $0$-cell $\rmod{A}$ is the identity functor
$id_{\rmod{A}}$. There is another bicategory which we denote by
$\mathscr{L}$ whose class of $0$-cells is the class of all rings
with local units, and with the same Home-Categories $\bara{ \rm
Funct}(\rmod{A},\rmod{B})$. Next, we formulate our results using the
bicategory $\mathscr{L}$ instead of $\td{\mathscr{L}}$.

\begin{proposition}\label{B-L}
Keep the above notations. There exists a bi-equivalence of
bicategories $\ring{F}: \ring{L} \to \ring{B}$ given locally by the
functors $$\xymatrix@R=0pt@C=50pt{ \ring{F}_{A,\,B}: \bara{ \rm
Funct}(\rmod{A},\rmod{B}) \ar@{>}[r] & \bimod{A}{B} \\ F \ar@{->}[r]
& F(A) \\ \eta_{-} \ar@{>}[r] & \eta_A, }$$ for every pair of ring
with local units $(A,B)$.
\end{proposition}
\begin{proof} It will be done in three steps.\\
\emph{Step.1. Homomorphism of bicategories}. The morphism $\ring{F}$
is the identity on the class of objects ($\ring{L}$ and $\ring{B}$
have the same class of objects). Taking two rings with local units
$A$ and $B$, the stated functor $\ring{F}_{A,\,B}$ is well defined
thanks to Lemmata \ref{ISO} and \ref{comonad1}(1). We need to show
its compatibility with the horizontal and vertical multiplications.
So let $C$ be another ring with local units and consider two
morphisms:  $\eta_{-}: F \to F'$ in the category $\bara{ \rm
Funct}(\rmod{A},\rmod{B})$, and $\zeta: G \to G'$ in the category
$\bara{ \rm Funct}(\rmod{B},\rmod{C})$. By Lemma \ref{comonad1}(1),
the morphism $$\xymatrix@C=60pt{ \ring{F}_{A,\,B}(F) .
\ring{F}_{B,\,C}(G) \,=\, F(A)\tensor{B} G(B)
\ar@{->}^-{\Theta_{F(A)}^G}[r] & \ring{F}_{A,\,C}(G\circ F)\,=\,
GF(A)}$$ is $A-C$--bilinear, where $\Theta_{F(A)}^G\,=\,
\left(\Upsilon_{F(A)}^G\right)^{-1}$ and $\Upsilon_{-}^-$ is the
natural transformation of Lemma \ref{ISO}. Moreover, applying Lemma
\ref{comonad1}, we obtain the following commutative diagram
$$\xymatrix@C=60pt@R=40pt{ F(A)\tensor{B}G(B)
\ar@{->}_-{\eta_A\tensor{B}\zeta_B\,=\,\ring{F}(\eta).\ring{F}(\zeta)}[d]
\ar@{->}^-{\Theta_{F(A)}^G}[r]  & GF(A)
\ar@{->}^-{\ring{F}(\eta.\zeta)\,=\, \zeta_{F'(A)} \circ G\eta_A}[d]
\\ F'(A)\tensor{B}G'(B) \ar@{->}^-{\Theta_{F'(A)}^{G'}}[r] & G'F'(A)
}$$ This implies that $\ring{V}_{F,G}:\,=\, \Theta_{F(A)}^G:
\ring{F}(F) . \ring{F}(G) \cong \ring{F}(G\circ F)$ is a natural
isomorphism. Clearly $\ring{F}(id_{\rmod{A}})\,=\,
id_{\rmod{A}}(A)\,=\, A$, and the compatibility with the
associativity, left and right multiplications by identities $1$-cell
is fulfilled.  Therefore, the pair $(\ring{F},\ring{V})$ establishes
a homomorphisms of bicategories from $\ring{L}$ to $\ring{B}$.

\emph{Step.2. Local equivalences of Hom-Categories}. Given two rings
with local units $A$ and $B$, and consider the stated functor
$\ring{F}_{A,\,B}$. Define the functor $\ring{G}_{A,\,B}:
\bimod{A}{B} \to \bara{ \rm Funct}(\rmod{B},\rmod{C})$ acting on
objects by $M \to -\tensor{A}M$ and on morphisms by $f \to
-\tensor{A}f$. It is clear that $\Upsilon_{-}^F$ gives a natural
isomorphism $ \ring{G}_{A,\,B} \circ \ring{F}_{A,\,B}(F)\,=\,
-\tensor{A}F(A) \cong F$, for any functor $F \in \bara{ \rm
Funct}(\rmod{A},\rmod{B})$. That is $\ring{G}_{A,\,B} \circ
\ring{F}_{A,\,B} \cong id_{\bara{ \rm Funct}(\rmod{A},\rmod{B})}$.
Conversely, for any $(A,B)$--bimodule $M$, we have a natural
isomorphism of bimodules $\ring{F}_{A,\,B} \circ \ring{G}_{A,\,B}(M)
\,=\, A\tensor{A}M \cong M$. That is $\ring{G}_{A,\,B} \circ
\ring{F}_{A,\,B} \cong id_{\bimod{A}{B}}$. Therefore,
$\ring{F}_{-,-}$ are equivalences of Hom-Categories.

\emph{Step.3. Surjectivity up to equivalences}. It is immediate since
$\ring{L}$ and $\ring{B}$ have the same class of objects which are
not altered by $\ring{F}$.
\end{proof}

\begin{remark}
As was pointed by the referee there is an alternative proof of
Proposition \ref{B-L} which uses results from
\cite{Grandjean/Vitale:1998}. There, firm rings (resp. firm modules)
(see Remark \ref{firm} for definitions) were termed regular algebras
(resp. regular modules). The regularity of a functor $F: \rmod{A}
\to \rmod{B}$ (\cite[Definition 1.5]{Grandjean/Vitale:1998}) between
categories of firm modules over firm rings $A$ and $B$, means that
$F(A)$ is a firm $A-B$-bimodule. Rings with local units are firm
rings. So by Lemma \ref{ISO}, right exact and direct sums preserving
functors (i.e., continuous) are regular. Hence Proposition \ref{B-L}
follows by \cite[Proposition 2.1]{Grandjean/Vitale:1998}. Another
way to obtain Proposition \ref{B-L} for rings with orthogonal
idempotents, is by using the arguments done before \cite[Section
2]{Grandjean/Vitale:1998}. For a given ring $A$ with a set of
orthogonal idempotents $\{e_i\}_{i \in I}$ (i.e., $A\,=\,\oplus_{i
\in I}Ae_i\,=\,\oplus_{i \in I}e_iA$), one can construct a right
$A$-linear map $\varphi: A \to \oplus_{a \in A}A$ such that
$$ \xymatrix{ A \ar@{->}^-{\varphi}[r] & \oplus_{a \in A}A
\ar@{->}^-{\theta}[r] & A\tensor{\mathbb{K}}A}$$ is a section for
the multiplication, where $\theta$ is the right $A$-linear map
induced be tensoring with a fixed element of $A$ (see \cite[page
142]{Grandjean/Vitale:1998}). Here $\varphi$ is defined as follows:
If $a \in A$, and $e_{i_1},\cdots,e_{i_n}$ the associated
idempotents such that $a \,=\, e_{i_1} a \,+\, \cdots\,+\,
e_{i_n}a$, we then take $\varphi(a) \,=\,e_{i_1} a \,\dotplus\,
\cdots\,\dotplus\, e_{i_n}a$.
\end{remark}

By $A\comonads{}$, we denote the category of all comonads in
$\rmod{A}$ whose underlying functors are right exact and preserve
direct sums, that is, the category of Definition \ref{continuos}
associated to the Grothendieck category $\rmod{A}$. The following
corollary is a direct consequence of Propositions \ref{B-L} and
\ref{comonad2}.
\begin{corollary}\label{corolario}
Let $A$ be a ring with local units. Then the functor
$$\ring{F}:\,\,A\comonads{}  \longrightarrow \coanillos{A}$$ defined by
$$\left(\underset{}{}\,
(F,\delta,\xi) \longrightarrow \left(\underset{}{}
F(A),\Upsilon_{F(A)}^F \,\circ\, \delta_A,\xi_A\right)\, \right),
\quad \left(\underset{}{}\, [\Phi: F \rightarrow F'] \longrightarrow
\left[\underset{}{} \Phi_A:F(A)\rightarrow F'(A) \right]\,\right) $$
establishes an equivalence of categories.
\end{corollary}

The right Eilenberg-Moore bicategory associated to $\ring{B}$ is
given by following corollary which is the non-unital version of
\cite[2.1]{Brzezinski/Kaoutit/Gomez:2004}

\begin{corollary}\label{bicatg} The following
data form a bicategory $\ring{R}$:
\begin{enumerate}
\item[$\bullet$]$0$-cells: corings $(\coring{C}:A)$ $($i.e., $A$
is a ring with local units and $\cC$ is an $A$--coring$)$.

\item[$\bullet$]$1$-cells: From $(\coring{D}:B)$ to
$(\coring{C}:A)$ are pairs $(M,\fk{m})$ consisting of an unital
$(A,B)$--bimodule $M$ and an $A-B$--bilinear map $\fk{m}:
\coring{C}\tensor{A}M \rightarrow M\tensor{B}\coring{D}$ compatible
with comultiplications and counits, that is, $\fk{m}$ satisfies
\begin{equation}\label{fks}
(M\tensor{B}\varepsilon_{\cD}) \,\circ\, \fk{m}  =
\varepsilon_{\coring{C}} \tensor{A} M, \,\, (\fk{m}\tensor{B}\cD)
\,\circ\, (\coring{C}\tensor{A}\fk{m}) \,\circ\,
(\Delta_{\coring{C}}\tensor{A}M) =  (M\tensor{B}\Delta_{\cD})
\,\circ\, \fk{m},
\end{equation}
where the first equality is up to the isomorphism $A\tensor{A}M
\cong M\tensor{B}B$. The identity $1$-cells for a given coring
$(\coring{C}:A)$ is given by the pair
$(A,\coring{C}\tensor{A}A \cong A\tensor{A}\coring{C})$.

\item[$\bullet$]$2$-cells: From $(M,\fk{m})$ to $(M',\fk{m}')$ are
$A-B$--bilinear maps $\fk{a}: \coring{C}\tensor{A}M \rightarrow M'$
satisfying
\begin{equation}\label{2-cells}
(\fk{a}\tensor{B}\coring{D}) \,\circ\, (\coring{C}\tensor{A}\fk{m})
\,\circ\, (\Delta_{\coring{C}} \tensor{A}M) \,\, = \,\, \fk{m}'
\,\circ\, (\coring{C}\tensor{A}\fk{a}) \,\circ\,
(\Delta_{\coring{C}} \tensor{A}M).
\end{equation}
\end{enumerate}
\end{corollary}

Laws composition are defined as in
\cite[2.1]{Brzezinski/Kaoutit/Gomez:2004}, and given by
\begin{equation}\label{tensor-1}
(M,\fk{m}) \boldsymbol{\tensor{}} (N,\fk{n}) \,\, =\,\,
\left(\underset{}{} M\tensor{B}N, \,(M\tensor{B}\fk{n}) \,\circ\,
(\fk{m}\tensor{B}N) \right)
\end{equation}
If $\fk{a}: \cC\tensor{A}M \rightarrow M'$ and $\fk{b}:
\cD\tensor{B}N \rightarrow N'$ are $2$-cells, then
\begin{equation}\label{tensor-2} \fk{a} \boldsymbol{\tensor{}} \fk{b} \,\, =\,\,
(M'\tensor{B}\fk{b}) \,\circ\, \fk{m'} \,\circ\,
(\cC\tensor{A}\fk{a}\tensor{B}N) \,\circ\,
(\Delta_{\cC}\tensor{A}M\tensor{B}N).
\end{equation}
\noindent The resulting category of all  $1-$ and $2$-cells from
$(\coring{D}:B)$ to $(\cC:A)$ is denoted by
${}_{(\coring{C}:\,A)}\ring{R}_{(\coring{D}:\,B)}$.

Let's keep now the notations of Section \ref{section0}. Then the right
Eilenberg-Moore $2$-category associated to $\ring{L}$ is defined as
follows

\begin{corollary}\label{2-cat}
The following data form a $2$-category $\ring{C}$:
\begin{enumerate}[$\bullet$]
\item $0$-\emph{cells}: They are pairs $(\bd{F}:A)$, that
is, $\bd{F}=(F,\delta,\xi) \in A\comonads{}$ where $A$ is a ring
with local units (i.e., $F$ is a continuous functor).

\item $1$-\emph{cells}: From $(\bd{G}:B)$ to
$(\bd{F}:A)$ is a pair $(S,\msf{s})$ consisting of a continuous
functor $S:\rmod{A} \rightarrow \rmod{B}$ and a natural transformation $\msf{s}: SF
\rightarrow GS$ satisfying the commutativity of diagrams in equation
\eqref{msf}

\item $2$-\emph{cells}: Given $(S,\msf{s})$ and $(S',\msf{s}')$
two $1$-cells from $(\bd{G}:B)$ to $(\bd{F}:A)$, a $2$-cell $\alpha:
(S,\msf{s}) \rightarrow (S',\msf{s}')$ is a natural transformation
$\alpha: SF \rightarrow S'$ satisfying the commutativity of diagram
in equation \eqref{dos-cells}.
\end{enumerate}
\end{corollary}
\noindent The category obtained by all $1-$ and $2$-cells from
$(\bd{G}:B)$ to $(\bd{F}:A)$ is denoted by
${}_{(\bd{F}:\,A)}\ring{C}_{(\bd{G}:\,B)}$.

The following is our main result of this section.

\begin{theorem}\label{F-G}
There is a bi-equivalence between the bicategory $\ring{R}$ of
corings over rings with locals units (Corollary \ref{bicatg}), and
the bicategory $\ring{C}$ whose objects are comonads with continuous
underlying functors over right unital modules (Corollary
\ref{2-cat}). This bi-equivalence is locally induced by the functors
$\bd{\ring{F}}_{(-,-)}$ defined by
$$\xymatrix@R=0pt@C=60pt{\bd{\ring{F}}_{(\bd{F},\,\bd{G})}:\,\,
{}_{(\bd{F}:\,A)}\ring{C}_{(\bd{G}:\,B)} \ar@{->}[r] &
{}_{(F(A):\,A)}\ring{R}_{(G(B):\,B)}  }$$
\begin{multline*}
(S, \msf{s}) \longrightarrow
 \left( \,S(A), \Upsilon_{S(A)}^{G} \,\circ\, \msf{s}_A
\,\circ\, \left( \Upsilon_{F(A)}^{S}\right)^{-1}
\, \right), \\
[ \alpha: SF \rightarrow S'] \longrightarrow \left[\, \alpha_A
\,\circ\, (\Upsilon_{F(A)}^{S}{})^{-1}: F(A)\tensor{A}S(A)
\rightarrow S'(A)\, \right],
\end{multline*}
where $\Upsilon_{-}^{-}$ are the natural isomorphisms of Lemma
\ref{ISO}, and where $(\bd{F}: A)$ (resp. $(\bd{G}:B)$) is sent to
$(F(A):A)$ (resp. to $(G(B):B)$) is the coring constructed in Proposition
\ref{comonad2}.
\end{theorem}
\begin{proof}
It is a consequence of Proposition \ref{B-L} and \cite[Remark
1.1]{Lack/Street:2002}.
\end{proof}

\begin{remark}
If we want to study any aspect in the right Eilenberg-Moore
bicategory $\ring{R}$, then it is convenient, using the local
equivalences $\bd{\ring{F}}_{-,-}$ stated in Theorem \ref{F-G}, to
transfer this study to the $2$-category $\ring{C}$. The local
equivalences in the other direction are given by the functors
$$\xymatrix@C=60pt{\bd{\ring{G}}_{(\coring{C},\,\coring{D})}:\,\,{}_{(\coring{C}:A)}\ring{R}_{(\coring{D}:B)}
\ar@{->}[r] & {}_{(\bd{F}:\,A)}\ring{C}_{(\bd{G}:\,B)} }$$ defined
by
\begin{multline*}
\left(\underset{}{} (M, \fk{m}) \longrightarrow
 (-\tensor{A}M, -\tensor{A}\fk{m}) \right), \\
\left(\underset{}{}[ \fk{a}: \coring{C}\tensor{A}M \rightarrow M']
\longrightarrow \left[\underset{}{}\, -\tensor{A}\fk{a}
:-\tensor{A}\coring{C}\tensor{A}M \rightarrow -\tensor{A}M'\,
\right] \right),
\end{multline*}
where $\bf{F}$ (resp. $\bf{G}$) is the comonad induced by the coring
$(\coring{C}:A)$ (resp. by $(\coring{D}:B)$). These local
equivalences are not in fact given by \cite[Remark
1.1]{Lack/Street:2002}. Their construction was given separately
using direct computations.
\end{remark}

\section{Base ring extension of a coring by an adjunction}\label{sect.app}
In this section we apply results from Sections \ref{section0} and
\ref{sectionI}, to extend the notion of base ring extension of a
coring by a (finitely generated and projective) module, introduced
in \cite{Brzezinski/Kaoutit/Gomez:2004}, to the case of rings with
local units. This will give a new class of corings over rings with
local units which includes some infinite comatrix corings
\cite{Kaoutit/Gomez:2004b}.
\bigskip

The following proposition characterizes an adjunction between right
unital modules with continuous functors (i.e. right exact and direct sums preserving functors).
\begin{proposition}\label{profinitos}
Let $A$ and $B$ be two rings with local units. The following
statements are equivalent.
\begin{enumerate}
\item[(i)] There is an adjunction $\xymatrix{S: \rmod{B}
\ar@<0,5ex>[r] & \ar@<0,5ex>[l] \rmod{A}:T}$ with $S$ left adjoint
to $T$, and such that $S$, $T$ are continuous functors.

\item[(ii)] There is an unital $(B,A)$--bimodule $\Sigma$ such
that $h\Sigma$ is finitely generated and projective unital right
$A$--module, for every $h \in \idemp{B}$.
\end{enumerate}
\end{proposition}
\begin{proof}
$(ii) \Rightarrow (i)$. We denote by
$\udual{\Sigma}=A\hom{A}{\Sigma}{A}B$ the unital right dual of
$\Sigma$. The bi-actions are defined as follows: For $a \in A$, $b
\in B$, and $\chi \in \Sigma^{\dag}$, we have $$ (a .\chi)(x)
\,\,=\,\, a\chi(x),\quad (\chi. b)(x)\,\,=\,\,\chi(bx),\qquad
\forall\, x \in \Sigma.$$ This is clearly the unital part of the
$(A,B)$--bimodule $\hom{A}{\Sigma}{A}$. The unit of the adjunction is given by
\begin{equation}\label{S-unit}
\xymatrix@R=0pt@C=60pt{ \eta_{Y_B}: Y \ar@{->}[r] &
Y\tensor{B}\Sigma\tensor{A}\Sigma^{\dag} \\ y \ar@{|->}[r] &
\sum_{i=1}^{n_h} y\tensor{B}u_i \tensor{A} (u_i^* \circ \pi_h) }
\end{equation}
where $h \in \idemp{B}$ such that $yh=y$, and
$\{(u_i,u^*_i)\}_{1\leq i \leq n_h} \subset h\Sigma \times
\hom{A}{h\Sigma}{A}$ is the finite right dual basis for $h\Sigma$,
and $\pi_h: \Sigma \to h\Sigma$ is the canonical projection. We
denote by $v_i^*$ the composition
$$\xymatrix@C=50pt{v_i^*: \Sigma \ar@{->}^-{\pi_h}[r] & h\Sigma
\ar@{->}^-{u_i^*}[r] & A } \,\, \in \udual{\Sigma}.$$ We claim that
$\eta_{Y_B} (yb)\,=\, \eta_{Y_B}(y)b$, for every $b \in B$ such that $b=hb=bh$. We have
\begin{eqnarray*}
  \sum_{i=1}^{n_h} yb\tensor{B}u_i\tensor{A}(u_i^*\circ \pi_h) &=&
  \sum_{i=1}^{n_h} y\tensor{B}bu_i\tensor{A}(u_i^*\circ \pi_h) \\
   &=& \sum_{i=1}^{n_h} \sum_{i'=1}^{n_h} y\tensor{B}u_{i'}\tensor{A} u_{i'}^*(bu_i)(u_i^*\circ \pi_h) \\
   &=& \sum_{i'=1}^{n_h} y\tensor{B}u_{i'}\tensor{A}\left(\sum_{i=1}^{n_h} u_{i'}^*(bu_i) (u_i^*\circ \pi_h)\right) \\
   &=& \sum_{i'=1}^{n_h} y\tensor{B}u_{i'}\tensor{A}
   \left(\left(\sum_{i=1}^{n_h} u_{i'}^*(bu_i) u_i^*\right)\circ \pi_h \right)\\
   &=& \sum_{i'=1}^{n_h} y\tensor{B}u_{i'}\tensor{A}((u_{i'}^*b) \circ \pi_h) \\
   &=& \sum_{i'=1}^{n_h} y\tensor{B}u_{i'}\tensor{A}(u_{i'}^* \circ
   \pi_h)b,\quad \text{ since } \,b=bh=hb,
\end{eqnarray*}
and this proves the claim. Next, we prove that $\eta_{Y_B}$ is
independent from the choice of $h$. So, let's fix an arbitrary
element $y\in Y$ and let $h' \in \idemp{B}$ be another unity for $y$
(i.e., $yh'=y=yh$). We consider as before $\{(x_j,w_j^*)\}_{1 \leq j
\leq n_{h'}} \subset h'\Sigma \times \Sigma^{\dag}$ the induced set
by the dual basis $\{(x_j,x_j^*)\}_{1 \leq j \leq n_{h'}}$ of
$h'\Sigma$, where $w_j^*= x_j^* \circ \pi_{h'}$. Henceforth, its
remains to prove that
$$\sum_{i=1}^{n_h} y\tensor{B}u_i \tensor{A} v_i^*  \,\, = \,\,
\sum_{j=1}^{n_{h'}} y\tensor{B}x_j \tensor{A} w_j^*. $$ So, let $h''
\in\idemp{B}$ such that $h=hh''=h''h$ and $h'=h'h''=h''h'$, and
consider again its corresponding set $\{(z_k,t_k^*)\}_{1 \leq k \leq
n_{h''}} \subset h''\Sigma \times \Sigma^{\dag}$, where $t_k^*=z_k^*
\circ \pi_{h''}$. Using elementary arguments, one can directly check
that
$$ v_i^* = (v_i^* \,\circ\, \tau_{h''}). h'' \text{ and } w_j^* =
(w_j^* \,\circ\, \tau_{h''}). h'', \text{ for every pair } (i,j) \in
\{1,\cdots,n_h\} \times \{1,\cdots,n_{h'}\},$$ where $\tau_{h''}:
h''\Sigma \rightarrow \Sigma$ is the canonical injection. On the
other hand, we have $\sum_{1 \leq k \leq n_{h''}} (v_i^* \circ
\tau_{h''})(z_k) t_i^*\,=\, (v^*_{i} \circ \tau_{h''}).h''$.  Taking
all these equalities into account, we compute
\begin{eqnarray*}
  \sum_{i}^{n_h} y\tensor{B}u_i \tensor{A} v_i^*
   &=& \sum_{i}^{n_h} y\tensor{B}u_i \tensor{A} (v_i^* \,\circ\, \tau_{h''})h''  \\
   &=& \sum_{i, k}^{n_h, n_{h''}} y\tensor{B}u_i \tensor{A} (v_i^* \,\circ\, \tau_{h''})(z_k) t_k^* \\
   &=& \sum_{k}^{ n_{h''}} y\tensor{B} \left( \sum_{i}^{n_h} u_i
   (v_i^* \,\circ\, \tau_{h''})(z_k)\right) \tensor{A} t_k^* \\
   &=& \sum_{k}^{ n_{h''}} y\tensor{B} \left( \sum_{i}^{n_h} u_i
   u_i^* (hz_k)\right) \tensor{A} t_k^*,\quad z_k \in h''\Sigma \\ &=& \sum_{k}^{ n_{h''}}
   y\tensor{B}hz_k\tensor{A}t_k^* \,\, = \,\, \sum_{k}^{ n_{h''}}
   y\tensor{B}z_k\tensor{A}t_k^*.
\end{eqnarray*}
Similar computation entails the equality $$
\sum_{j}^{n_{h'}}y\tensor{B}x_j \tensor{A} w_j^* \,\, =\,\,
\sum_{k}^{ n_{h''}} y\tensor{B}z_k\tensor{A}t_k^*,$$ and this proves
the desired independence. Therefore, $\eta_{Y_B}$ is a well defined
right $B$-linear map, for every right unital $B$--module $Y$.
Clearly, $\eta_{-}: id_{\rmod{B}}(-) \rightarrow
-\tensor{B}\Sigma\tensor{A}\Sigma^{\dag}$ is a natural
transformation. The counit is given by
\begin{equation}\label{S-counit}
\xymatrix@R=0pt@C=60pt{ \zeta_{X_A}:
X\tensor{A}\Sigma^{\dag}\tensor{B}\Sigma \ar@{->}[r] & X \\
x\tensor{A}\varphi\tensor{B}u \ar@{|->}[r] & x\varphi(u).}
\end{equation}
Lastly, one can easily show that $$ (\zeta_{X}
\tensor{A}\Sigma^{\dag}) \,\circ\, \eta_{X\tensor{A}\Sigma^{\dag}} =
X\tensor{A}\Sigma^{\dag},\text{ and }\,\, \zeta_{Y\tensor{B}\Sigma}
\,\circ\, (\eta_{Y} \tensor{B}\Sigma) = Y \tensor{B}\Sigma,$$ for
every pair of unital modules $(Y_B,X_A)$ which implies the desired
adjunction taking $S(-)=-\tensor{B}\Sigma$ and
$T(-)=-\tensor{A}\Sigma^{\dag}$.

$(i) \Rightarrow (ii)$. By Lemma \ref{ISO} we know that ${}_BS(B)_A$
and ${}_AT(A)_B$ are unital bimodules, and that $S(-) \cong
-\tensor{B}S(B)$, $T(-) \cong -\tensor{A}T(A)$ are natural
isomorphisms. Taking ${}_B\Sigma_A = S(B)$, and $h \in \idemp{B}$, we deduce that
$h\Sigma = hS(B) \cong hB\tensor{B}S(B) \cong S(hB)$ is a right
$A$--linear isomorphism. Henceforth, it remains to show that
$S(hB)$ is a finitely generated and projective module, for every $h\in
\idemp{B}$. So, the natural isomorphism of the stated adjunction
gives us the following chain of natural isomorphisms
$$\hom{A}{S(hB)}{-} \,\,\cong \,\, \hom{B}{hB}{T(-)} \cong T(-)h
\cong -\tensor{A}T(A)h, $$ for every $h\in \idemp{B}$. That is, the
functor $\hom{A}{S(hB)}{-}$ preserves inductive limits, and so
$S(hB)$ is a finitely generated and projective $A$--module for every
$h \in \idemp{B}$.
\end{proof}

\begin{remark}
Considering $\Sigma$ an unital $(B,A)$--bimodule, we can easily check,
using the partial ordering on idempotent elements, that
$\dlimit{{}_h(h\Sigma)}\,\cong\, \Sigma$ as right unital
$A$-modules. In fact $\{h\Sigma\}_{h\in \,\idemp{B}}$ is a split
direct system of right unital $A$-module (see \cite[Section
1]{Vercruysse:2004}). If we assume that $\Sigma$ satisfies condition
$(ii)$ of Proposition \ref{profinitos}, then $\Sigma_A$ is locally
projective in the sense of \cite{Anh/Marki:1987},
equivalently, it is strongly locally projective in the sense of
\cite[Theorem 2.17]{Vercruysse:2004}.
\end{remark}

\begin{remark}\label{tenso-dos}
Let $A$, $B$, and $C$ be rings with local units. Consider $\Sigma$
(respectively $W$) an unital $(B,A)$--bimodule (respectively
$(C,B)$--bimodule) such that $h \Sigma$ (respectively $gW$) is
finitely generated and projective unital right $A$--module
(respectively $B$--module), for every $h \in \idemp{B}$
(respectively $g \in \idemp{C}$). Then $W\tensor{B}\Sigma$ is an
unital $(C,A)$--bimodule such that $g(W\tensor{B}\Sigma)$ is
finitely generated and projective unital right $A$--module, for
every $g \in \idemp{C}$. Furthermore, if we put
$\Sigma^{\dag}=A\hom{A}{\Sigma}{A}B$ (respectively
$W^{\dag}=B\hom{B}{W}{B}C$), then
\begin{equation}\label{ismo-prof}
(W\tensor{B}\Sigma)^{\dag}=A\hom{A}{W\tensor{B}\Sigma}{A}C \,\,
\cong \,\, \Sigma^{\dag} \tensor{B} W^{\dag}
\end{equation}
is an isomorphism of unital $(A,C)$--bimodules. Effectively, let $g
\in \idemp{C}$ any idempotent element, so there exists an unital
right $B$--module $N$ such that
$$ gW \,\oplus \, N = \bigoplus_{i=1}^{n}h_iB$$ where each $h_i \in
\idemp{B}$. Applying the tensor product $-\tensor{B}\Sigma$, we
obtain $$ (gW\tensor{B}\Sigma)\,\oplus\, (N\tensor{B}\Sigma) \cong
\bigoplus_{i=1}^n(h_i\Sigma)$$ an isomorphism of unital right
$A$--modules. Since the right hand module is a finitely generated and
projective $A$--module, we get that $gW\tensor{B}\Sigma$ is also
finitely generated and projective as an $A$--module, and this proves the
first claim. Now, using the adjunctions arising from the proof of
Proposition \ref{profinitos} and the usual Hom-Tensor adjunction,
we get the isomorphism of equation \eqref{ismo-prof}.
\end{remark}

It is convenient to adopt the notations of the proof of
Proposition \ref{profinitos}. Thus, if $\Sigma$ is any $(B,A)$--bimodule
we denote by $\Sigma^{\dag}=A\hom{A}{\Sigma}{A}B$. When
$h\Sigma$ is a finitely generated and projective right $A$--module,
for some $h \in \idemp{B}$, we consider the set $\{(u_i,v_i^*)\}_{1
\leq i \leq n_h} \subset h\Sigma \times \Sigma^{\dag}$ where
$\{(u_i,u_i^*)\}_{1 \leq i \leq n_h} \subset h\Sigma \times
\hom{A}{h\Sigma}{A}$ is the finite dual basis for $h\Sigma$, where
$v_i^*= u_i^* \,\circ\, \pi_h$ and $\pi_h: \Sigma \rightarrow
h\Sigma$ is the canonical projection.

\begin{corollary}\label{otros}
Let $A$ and $B$ be two rings with local units together with an
unital $(B,A)$--bimodule $\Sigma$ and a $B$--coring
$(\coring{D},\Delta_{\coring{D}},\varepsilon_{\coring{D}})$. Assume
that $h\Sigma$ is finitely generated and projective module for every
$h \in \idemp{B}$. Then the unital $A$--bimodule $\Sigma^{\dag}
\tensor{B} \coring{D} \tensor{B} \Sigma$ admits the structure of an
$A$--coring with comultiplication defined by
\begin{multline*}
\Delta: \Sigma^{\dag} \tensor{B} \coring{D} \tensor{B} \Sigma
\longrightarrow \Sigma^{\dag} \tensor{B} \coring{D} \tensor{B}
\Sigma \tensor{A} \Sigma^{\dag} \tensor{B} \coring{D} \tensor{B}
\Sigma \\ \varphi\tensor{B}d\tensor{B}u \longmapsto \sum_{i,(d)}
\varphi \tensor{B} d_{(1)} \tensor{B} u_i \tensor{A} v_i^*
\tensor{B} d_{(2)} \tensor{B} u,
\end{multline*}
where $\{(u_i,v_i^*)\}_{i} \subset h\Sigma \times \Sigma^{\dag}$ is
the finite set induced by the dual basis of $h\Sigma$, where $h \in
\idemp{B}$ such that $\varphi h= \varphi$, $hu=u$, $d=hd=dh$, and
where $\Delta_{\coring{D}}(d) = \sum_{(d)}d_{(1)}
\tensor{B}d_{(2)}$. The counit is defined by
$$
\varepsilon: \Sigma^{\dag} \tensor{B} \coring{D} \tensor{B} \Sigma
\longrightarrow A,\quad \left(\varphi\tensor{B}d \tensor{B} u
\longmapsto \varphi(\varepsilon_{\coring{D}}(d)u)\right). $$
\end{corollary}
\begin{proof}
By Proposition \ref{profinitos}, we know that $-\tensor{B}\Sigma:
\rmod{B} \rightarrow \rmod{A}$ is left adjoint to
$-\tensor{A}\Sigma^{\dag}: \rmod{A} \rightarrow \rmod{B}$ with unit
$\eta_{-}: id_{\rmod{B}}(-) \rightarrow
-\tensor{B}\Sigma\tensor{A}\Sigma^{\dag}$ and counit $\zeta_{-}:
-\tensor{A}\Sigma^{\dag}\tensor{B}\Sigma \rightarrow
id_{\rmod{A}}(-)$ given explicitly by equations \eqref{S-unit} and
\eqref{S-counit}. Applying Lemma \ref{dos}(1) to the comonad
$(-\tensor{B}\coring{D}, -\tensor{B}\Delta_{\coring{D}},
-\tensor{B}\varepsilon_{\coring{D}})$ in $\rmod{B}$, we obtain a new
comonad $(G,\Omega,\gamma)$ in $\rmod{A}$ with
\begin{eqnarray*}
  G &=& -\tensor{A}\Sigma^{\dag}\tensor{B}\coring{D}\tensor{B}\Sigma \\
  \Omega_{-} &=& \left( \underset{}{}\eta_{-\tensor{A}\Sigma^{\dag}\tensor{B}\coring{D}}
\tensor{B}\coring{D} \tensor{B} \Sigma \right)\,\circ\,
\left(\underset{}{}-\tensor{A}\Sigma^{\dag}\tensor{B}\Delta_{\coring{D}}\tensor{B}\Sigma\right) \\
  \gamma_{-} &=& \zeta_{-}\,\circ\,
\left(\underset{}{}-\tensor{A}\Sigma^{\dag}\tensor{B}\varepsilon_{\coring{D}}
\tensor{B}\Sigma\right)
\end{eqnarray*}
Now, Proposition \ref{comonad2} implies that $G(A)$ admits the
structure of an $A$--coring. Since $G(A) \cong \Sigma^{\dag}
\tensor{B}\coring{D} \tensor{B} \Sigma$ is obviously an isomorphism
of unital $A$--bimodules, this structure can be transferred to
$\Sigma^{\dag} \tensor{B}\coring{D} \tensor{B} \Sigma$ with
comultiplication and counit computed explicitly from the maps
$\Omega_A$ and $\gamma_A$. This, in fact, leads exactly to the stated
structure.
\end{proof}

Recently, in  \cite{Gomez/Vercruysse:2005} and
\cite{Caenepeel/DeGroot/Vercruysse:2005} new generalizations of
infinite comatrix corings, earlier introduced in
\cite{Kaoutit/Gomez:2004b}, were given in the context of non-unital
(firm) rings. The following example gives another way to construct infinite comatrix corings by
using Corollary \ref{otros}.

\begin{example}
Assume that $A$ is a ring with identity $1_A$, and denote by
$add(A_A)$ the full subcategory of all finitely generated and
projective unital right $A$--modules. Consider a
$\Bk$--additive small category $\cat{A}$ and its induced ring with enough
orthogonal idempotents $B = \oplus_{(\fk{p},\fk{q}) \in \cat{A}^2}
\hom{\cat{A}}{\fk{p}}{\fk{q}}$: These are $\{1_{\fk{p}}\}_{\fk{p}
\in \cat{A}} \subset B$, where each of the $1_{\fk{p}}$'s is the
image of the identity $1_{\mathrm{End}_{\cat{A}}(\fk{p})}$. Given an
additive faithful functor $\omega: \cat{A} \rightarrow add(A_A)$, we
get an unital $(B,A)$--bimodule $\Sigma=\oplus_{\fk{p} \in \cat{A}}
\omega(\fk{p})$. It is clear that $1_{\fk{p}}\Sigma =
\omega(\fk{p})$, for every object $\fk{p} \in \cat{A}$. Therefore,
$h\Sigma$ is a finitely generated and projective right $A$--module,
for every $h \in \idemp{B}$. Finally, considering $B$ as a trivial
$B$--coring, we obtain by Corollary \ref{otros} an $A$-coring
$\Sigma^{\dag} \tensor{B}B \tensor{B}\Sigma \cong \Sigma^{\dag}
\tensor{B}\Sigma$, where $\Sigma^{\dag}=\oplus_{\fk{p}
\in \cat{A}} \hom{A}{\omega(\fk{p})}{A}$.
\end{example}

\begin{remark}\label{firm}
Rings with local units are in fact a sub-class of firm rings. Recall
that an associative ring $R$ is firm if the multiplication
$R\tensor{R}R \to R$ is an isomorphism. Unital modules are extended
to firm modules, i.e. a right $R$--module $M$ with action
$M\tensor{R}R \to M$ an isomorphism of right $R$-modules. The
results of this paper can be extended to this class of rings by
using the categorical version of Lemmata \ref{comonad1} and
\ref{ISO} stated in \cite[39.3, 39.5]{Brzezinski/Wisbauer:2003} with
firm base rings. The fact that right exact and direct sums
preserving functors between firm modules are naturally isomorphism
to the tensor product functors, has been recently proved by J.
Vercruysse in \cite[Theorem 3.1]{Vercruysse:2006} (see also
\cite[Proposition 1.6]{Grandjean/Vitale:1998}). A characterization
of an adjunction whose both functors are right exact and preserve
direct sums (as in Proposition \ref{profinitos}) was extended to the
case of firm modules in \cite[Theorem 2.4]{Vercruysse:2006}.
\end{remark}

\section*{Acknowledgements}
I would like to thank my supervisor Professor J. G\'omez-Torrecillas
for helpful comments. I also would like to thank the referee for her
or his suggestions, for pointing out a short proof of Theorem
\ref{F-G}, and also for drawing my attention to the reference
\cite{Grandjean/Vitale:1998}. This research was supported by grants
MTM2004-01406 and P06-FQM-1889.


\begin{thebibliography}{10}
\bibitem{Abrams:1983}
G.~D. Abrams, \emph{Morita equivalence for rings with local units},
Commun.  Algebra \textbf{11} (1983), no.~8, 801--837.

\bibitem{Anh/Marki:1983}
P.N. \'{A}nh and L.~{M}\'{a}rki, \emph{Rees matrix rings}, J.
Algebra \textbf{81} (1983), 340--369.

\bibitem{Anh/Marki:1987}
P.N. \'{A}nh and L.~{M}\'{a}rki,  \emph{Morita equivalence for rings
without identity}, Tsukuba J. Math.  \textbf{11} (1987), no.~1,
1--16.

\bibitem{Bautista/Colavita/Salmeron:1981}
R.~Bautista, L.~Colavita, and L.~Salmer\'{o}n, \emph{On adjoint
functors in representation theory}, Lecture Notes in Math., vol.
903, Springer-Verlag, 1981, pp.~9--25.

\bibitem{Beck:1969}
J.~Beck, \emph{Distributive laws}: In Seminar on Triples and
Categorical Homology Theory. Lecture Notes in Math., vol. 80,
Springer-Verlag, 1969, pp.~119--140.

\bibitem{Benabou:1967}
J.~B\'enabou, \emph{Introduction to bicategories}: In Report of the
{M}idwest {C}ategory {S}eminar. Lecture Notes in Math., vol. 47,
Springer-Verlag, 1967, pp.~1--77.

\bibitem{Brzezinski/Kaoutit/Gomez:2004}
T.~Brzezi\'{n}ski, L.~El~Kaoutit, and J.~G\'{o}mez-Torrecillas,
\emph{The bicategories of corings}, J. Pure Appl. Algebra
\textbf{205} (2006), 510--541.

\bibitem{Brzezinski/Wisbauer:2003}
T.~Brzezi\'nski and R.~Wisbauer, Corings and Comodules, \emph{London
Math. Soc. Lect. Note Ser.} \textbf{309}, Cambridge University
Press, Cambridge, 2003.

\bibitem{Caenepeel/DeGroot/Vercruysse:2005}
S.~Caenepeel, E.~De Groot, and J.~Vercruysse, \emph{Constructing
infinite comatrix corings from colimits}, Appl. Categor. Struct.
\textbf{14} (2006), 539--565.

\bibitem{Iglesias/Gomez/Nastasescu:1999}
F.~Casta\~{n}o Iglesias, J.~G\'{o}mez-Torrecillas, and
C.~$\mathrm{N}${\v{a}}st\v{a}sescu, \emph{{F}robenius functors:
Aplications}, Commun. Algebra \textbf{27} (1999), no.~10,
4879--4900.

\bibitem{Eilenberg/Moore:1965}
S.~Eilenberg and J.~C. Moore, \emph{Adjoint functors and triples},
Illinois J. Math. \textbf{9} (1965), 381--398.

\bibitem{Kaoutit/Gomez:2003a}
L.~El~Kaoutit and J.~G\'{o}mez-Torrecillas, \emph{Comatrix corings:
{G}alois coring, {D}escent theory, and a structure theorem for
cosemisimple corings}, Math. Z. \textbf{244} (2003), 887--906.

\bibitem{Kaoutit/Gomez:2004b}
L.~El~Kaoutit and J.~G\'{o}mez-Torrecillas,, \emph{Infinite comatrix
corings}, Int. Math. Res. Notices \textbf{39} (2004), 2017--2037.

\bibitem{Kaoutit/Gomez/Lobillo:2004c}
L.~El~Kaoutit, J.~G\'{o}mez-Torrecillas, and F.~J. Lobillo,
\emph{Semisimple corings}, Algebra Colloq. \textbf{11} (2004),
no.~4, 427--442.

\bibitem{Fisher/Newell:1971}
J.~Fisher ~Palmquist and D.~C.~ Newell, \emph{Bifunctors and adjoint
pairs}, Trans. Amer. Math. Soc. \textbf{155} (1971), no.~2,
293--303.

\bibitem{Fuller/Hullinger:1978}
K.~R.~ Fuller and H~Hullinger, \emph{Rings with finiteness
conditions and their categories of functors}, J. Algebra \textbf{55}
(1978), 94--105.

\bibitem{Gabriel:1962}
P.~Gabriel, \emph{Des cat\'{e}gories ab\'{e}liennes}, Bull. Soc.
Math. France \textbf{90} (1962), 323--448.

\bibitem{Gabriel/Popescu:1964}
P.~Gabriel and N.~Popescu, \emph{Caract\'{e}risation des
cat\'{e}gories ab\'{e}lienne avec g\'{e}n\'{e}rateurs et limites
inductives exactes}, C. R. Acad. Sci. Paris \textbf{258} (1964),
4188--4191.

\bibitem{Gomez/Vercruysse:2005}
J.~G\'{o}mez-Torrecillas and J.~Vercruysse, \emph{Comatrix corings and
{G}alois comodules over firm rings}, Algebra Rep. Theory \textbf{10} (2007), 271--306.

\bibitem{Grandjean/Vitale:1998}
F.~Grandjean and E.M.~Vitale, \emph{Morita equivalence for regular
algebras}, Cahiers Topologie G\'{e}om. Diff\'{e}rentielle Cat\'{e}g.
\textbf{XXXIX} (1998), no. 2, 137--153.

\bibitem{Guzman:1989}
F.~Guzman, \emph{Cointegration, {R}elative {C}ohomology for
{C}omdules, and  {C}oseparable {C}oring}, J. Algebra \textbf{126}
(1989), 211--224.

\bibitem{Huber:1961}
P.~J. Huber, \emph{Homotopy theory in general categories}, Math.
Ann.  \textbf{144} (1961), 361--385.

\bibitem{Kleisli:1965}
H.~Kleisli, \emph{Every standard construction is induced by a pair
of adjoint functors}, Proc. Amer. Math. Soc. \textbf{16} (1965),
544--546.

\bibitem{Lack/Street:2002}
S.~Lack and R.~Street, \emph{The formal theory of monads {II}}, J.
Pure Appl.  Algebra \textbf{175} (2002), 243--265.

\bibitem{Mitchell:1965}
B.~Mitchell, \emph{Theory of {C}ategories}, {N}ew {Y}ork and
{L}ondon ed., Academic {P}ress, 1965.

\bibitem{Popescu:1973}
N.~Popescu, \emph{Abelian {C}ategory with {A}pplication to {R}ing
and {M}odules}, Academic {P}ress {L}ondon and {N}ew {Y}ork, 1973.

\bibitem{Street:1972}
R.~Street, \emph{The formal theory of monads}, J. Pure Appl. Algebra
\textbf{2}  (1972), 149--168.

\bibitem{Sweedler:1975}
M.~Sweedler, \emph{The predual theorem to the {J}acobson-{B}ourbaki
theorem}, Tran. Amer. Math. Soc. \textbf{213} (1975), 391--406.

\bibitem{Vercruysse:2004}
J.~ Vercruysse, \emph{Local units versus local projectivity.
{D}ualisations: {C}orings with local structure maps}, Commun.
Algebra \textbf{34} (2006), 2683--2711.

\bibitem{Vercruysse:2006}
J.~ Vercruysse, \emph{Equivalences between categories of modules and
categories of comodules}, arXiv:math.RA/0604423.

\bibitem{Watts:1960}
C.~E. Watts, \emph{Intrinsic characterization of some additive
functors}, Proc. Amer. Math. Soc. \textbf{11} (1960), 5--8.

\end{thebibliography}
\end{document}